# A Dynamic Likelihood Approach to Filtering for Advection-Diffusion Dynamics

Johannes Krotz,[a] Juan M. Restrepo,[a,b] and Jorge Ramirez,[b]

[a] *Department of Mathematics, University of Tennessee, Knoxville, Knoxville TN, USA*
[b] *Computer Science and Mathematics Division, Oak Ridge National Laboratory, Knoxville TN 37831, USA*

ABSTRACT: A Bayesian data assimilation scheme is formulated for advection-dominated advective and diffusive evolutionary problems, based upon the Dynamic Likelihood (DLF) approach to filtering. The DLF was developed specifically for hyperbolic problems –waves–, and in this paper, it is extended via a split step formulation, to handle advection-diffusion problems. In the dynamic likelihood approach, observations and their statistics are used to propagate probabilities along characteristics, evolving the likelihood in time. The estimate posterior thus inherits phase information. For advection-diffusion the advective part of the time evolution is handled on the basis of observations alone, while the diffusive part is informed through the model as well as observations. We expect, and indeed show here, that in advection-dominated problems, the DLF approach produces better estimates than other assimilation approaches, particularly when the observations are sparse and have low uncertainty. The added computational expense of the method is cubic in the total number of observations over time, which is on the same order of magnitude as a standard Kalman filter and can be mitigated by bounding the number of forward propagated observations, discarding the least informative data.

## 1. Introduction

A general framework in Bayesian estimation to assimilate observations and model predictions has become known as data assimilation. Models are used to inform a prior and observations inform the likelihood. For time-dependent problems, the estimation objective is to find the evolution of moments of the posterior of a time-dependent state variable, conditioned on observations. A variety of computational methodologies have been proposed to accomplish this (see Sarkka (2013) and references therein). In linear problems with Gaussian noise processes, the variance minimizer estimate of the time-dependent mean and variance of the posterior can be obtained sequentially by the Kalman Smoother Rauch et al. (1965) or partially by the Kalman Filter (KF) Kalman (1960). Kushner, Stratanovich, Pardoux (see, for example, Kushner (1962)) proposed a variance minimizer estimate for the nonlinear/non-Gaussian problem, however, it is computationally tractable only for very low-dimensional state variable problems. Successful approximations of the estimate can sometimes be obtained via generalizations like the Extended Kalman Filter McElhoe (1966); Smith et al. (1962), or the Unscented Kalman Filter Julier and Uhlmann (2004), among others. Sample estimates can be approximated via the Ensemble Kalman Filter Evensen (1994, 2004) and its variants, the path integral method Alexander et al. (2005), and various particle filter schemes Moral et al. (2006); Chorin et al. (2010); Restrepo and Ramírez (2022). There are estimators that have special properties (see Rosenthal et al. (2017)) or that exploit the underlying dynamics of the problem. An example of the latter is the the dynamic likelihood filtering approach (DLF), first proposed in Restrepo (2013).

*Corresponding author*: Juan M. Restrepo, restrepojm@ornl.gov

The DLF is denoted an "approach" rather than a filtering method because, in principle, it applies to any of the linear or nonlinear data assimilation methodologies. It was developed specifically for problems in wave dynamics, in general, hyperbolic partial differential equations. The crux of the DLF approach is to modify the conditional, posterior distribution of the state variable by exploiting a property peculiar to wave problems: finite-time propagation of information, which is utilised to propose a dynamic likelihood. A second aspect of DLF is that it tracks the state variables of the partial differential equation along characteristics, thus obtaining stochastic differential equations. Peculiarities of the DLF approach are that phase information enters directly into the estimation, and that we can make Bayesian estimates at times when observations are available and when they are not (even in the near future). The main practical advantage of the method is that it addresses the more common situation in wave problems: sparse observation networks that are, nevertheless, fairly low in noise. Under these circumstances, as was shown in Foster and Restrepo (2022), the DLF produces superior estimates when compared to the best traditional estimator.

In this paper, we develop the dynamic likelihood approach for data assimilation problems in transport modeled by forced advection-diffusion equations. We thus expand the range of applicability of this estimation approach to an important class of dynamics. We will focus on finding estimates of quantities of interest when the source of uncertainties appears in the advection process and the forcing. The statement of the problem appears in Section 2. The DLF evolves the likelihood forward in time along characteristics by generating *pseudo-observations* at times between actual observations. A pseudo-observation is derived from a real observation at a previous time. The pseudo-observation framework appears in Section a. This section also details how the DLF approach applies to the advection-diffusion dynamics, using techniques and ideas similar to a Kalman Filter.

To appreciate the practicality of the methodology, we present in Section 4 an accounting of the cost of implementing the DLF on a sequential estimator data assimilation method. In Section b, we compare the DLF approach proposed for advection-diffusion problems to the outcomes obtained via a Kalman Filter because the Kalman estimates for this problem are familiar, optimal, and easily understood. A discussion and conclusions appear in Section 6.

## 2. Statement of the Problem

At issue is the estimation of the posterior covariance of a noisy scalar state variable $u(x,t)$ given noisy observations, and the minimization of its trace. Here and throughout, $x$ denotes space and $t$, time, The state variable obeys a noisy advection-diffusion initial value problem. We will develop a DLF approach to a particular filtering estimation scheme. We focus on linear dynamics since it allows us to evaluate the DLF approach in comparison with optimal filtering schemes, nevertheless, we argue that the development presented herein will extend to a number of nonlinear cases.

We are motivated to consider the DLF approach to the dynamics of advection-diffusion because it was shown in Foster and Restrepo (2022) that for hyperbolic dynamics, the DLF returned significantly better estimates on noisy hyperbolic problems, particularly when the observations were sparse yet had low uncertainty –which is the more common practical situation. We will, in fact, show that for the advection-dominated case, the DLF approach yields better estimates than other filtering approaches.

Since we are specializing to the linear advection-diffusion initial value problem with known Gaussian noise processes it is possible to fully determine the posterior distribution with the determination of the posterior mean and variance. We connote a sample time series from the distribution of $u(x,t)$ as the **truth**. We will make use of the truth for testing the performance of the DLF. In practice, the truth is not available to us. Instead, we are given an approximate solution of the stochastic advection-diffusion initial value problem, with known errors, often in the form of a computer code. The estimation problem will thus be one of finding moments of the posterior model state variable, given observations.

### *a. Dynamics, Model, Observations*

The space interval over which the dynamical system is defined will be $[0,L] \subset \mathbb{R}$, with periodic boundaries. Space will be discretized by a grid $X$ of equidistant nodes $X = \left\{x^k = k \cdot \Delta x\right\}_{k=0}^{K}$, with $x^K = L - \Delta x$ due to the periodic boundary conditions. The time interval shall be $[0,t_N] \subset \mathbb{R}$ discretized in equal time steps $T := \{t_n = n\Delta t\}_{n=0}^{N}$. We will denote the set $T$ as *estimation times*. On $[0,L] \times [0,t_N]$ we consider the random field $u(x,t)$, which obeys the stochastic initial value problem

$$u_t - C(x,t)u_x = Du_{xx} + F(x,t), \quad t > 0, \quad x \in [0,L], \\ u(x,0) = u_0(x), \quad x \in [0,L]. \tag{1}$$

The subscripts $x$ and $t$ connote partial differentiation with respect to these variables. The periodic initial condition is $u_0(x)$ is known or is drawn from an assumed known probability distribution $\mathcal{U}$. The parameter $D$ is the diffusion constant, $F(x,t)$ and $C(x,t)$ are forcing and wave speed terms, respectively. It will be assumed that

$$C(x,t) = c(x,t) + \phi(x,t), \tag{2}$$

$$F(x,t) = f(x,t) + \chi(x,t), \tag{3}$$

where $f(x,t)$ and $c(x,t)$ are the forcing and phase speed, respectively, which are assumed known, deterministic and periodic on $[0,L]$. The random fields $\phi$ and $\chi$ have the form $\phi(x,t)dt = AdW^c(x,t)$, $\chi(x,t)dt = BdW^u(x,t)$ with $A$ and $B$ known constants, $dW^c$ and $dW^u$ incremental zero-mean Wiener processes, assumed uncorrelated. We define a semi-continuous ensemble member solution to (1) on the space grid as $\boldsymbol{U}(t) := \left(U^k(t)\right)_{k=0}^{K} = \left(u(X,t)\right)_{k=0}^{K}$.

Going forward, bold variables will denote vectors or matrices. Superindices are space, subindices are time. For all variables with a single index, e.g. $a_i$, we denote by $a_{k:n} = \{a_i\}_{i=k}^{n}$ the union over all indices between $k$ and $n$.

We connote $v$ as the **model** approximation to (1). We will build a specific one here as follows: on the grid $X \times T$, the values of $v$ are obtained by a forward numerical solution of the SDE (1). At each time $t_n$, we denote the collection of values of the model $v(\cdot,t_n)$ on $X$ by

$$\boldsymbol{V}_n = v(X,t_n) \tag{4}$$

The vector $\boldsymbol{V}_n$ is evolved forward with the SDE solver

$$\boldsymbol{V}_{n+1} = \boldsymbol{L}_n \boldsymbol{V}_n + \sqrt{\Delta t}\,\Delta \boldsymbol{w}_n + \Delta t\, \boldsymbol{f}_n \tag{5}$$

where $\boldsymbol{L}_n \in \mathbb{R}^{K \times K}$ is a numerical operator approximating the linear terms in (1), $\boldsymbol{f}_n := \left(f(x^k,t_n)\right)_{k=0}^{K}$ and $\Delta \boldsymbol{w}_n \in \mathbb{R}^N$ is mean-zero Gaussian vector with covariance matrix $\boldsymbol{Q}_n$ accounting for the stochasticity of (1) and the model error. The distribution of $\boldsymbol{V}_0$ and $\langle \boldsymbol{w}_n \boldsymbol{w}_n{}^\top \rangle$ are assumed known. For $x$-values that are off-grid, we use *linear interpolation in space* to extend the outputs of the numerical SDE solver (5) to $[0,L] \times T$. Namely, for $x \notin X$, we define

$$v(x,t_n) = \boldsymbol{H}(x)\boldsymbol{V}_n, \quad n = 1,\ldots,N \tag{6}$$

where $\boldsymbol{H}$ is a linear interpolation operator to be specified later. The model, up to time $t_{n_m}$, is used to inform the prior $\pi(\boldsymbol{V}_{0:n_m})$.

In practice, observations are obtained from instruments and the error is instrument-related. Here we generate them synthetically from the truth. We assume that **observations** are available at observation times $\{t_{n_1}, \ldots, t_{n_M}\} =: T_O \subset T$. The set of available observations is $O := \{(\boldsymbol{y}_m, \boldsymbol{Y}_m)\}_{m=1}^M$ which provide, up to noise, temporally and spatially localized records on the value of $u$. Specifically, the observation pair $(\boldsymbol{y}_m, \boldsymbol{Y}_m)$ corresponds to a time $t_{n_m}$ in the set of observation times $T_O \subset T$. The vector $\boldsymbol{y}_m \in [0, L]^I$ contains the $I \in \mathbb{N}$ locations where the observations were recorded, and $\boldsymbol{Y}_m \in \mathbb{R}^I$ is a measurement of the value of $u$ at those locations at time $t_{n_m}$. Specifically,

$$Y_m^i = u(y_m^i, t_{n_m}) + \epsilon_m^i, \quad i = 1, \ldots, I \tag{7}$$

where the measurement error $\boldsymbol{\epsilon}_m$ is a mean-zero, normal vector in $\mathbb{R}^I$ with known covariance. Note that observation times $T_O \subset T$ do not include all times in $T$ and that the number of observations $I$ does not depend on $t$.

By Bayes Law

$$\pi(\boldsymbol{V}_{0:n_m}|\boldsymbol{Y}_{1:m}) \propto \pi(\boldsymbol{Y}_{1:m}|\boldsymbol{V}_{0:n_m})\pi(\boldsymbol{V}_{0:n_m}). \tag{8}$$

The likelihood at time $t_{n_m} \in T_O$, informed by observations is $\pi(\boldsymbol{Y}_{1:m}|\boldsymbol{V}_{0:n_m})$. The prior $\pi(\boldsymbol{V}_{0:n_m})$ is informed by the model.

### *b. The Kalman Filter (KF)*

In principle, the DLF approach applies to most sequential filtering schemes. Since the problem we are considering is linear, we will be testing the DLF approach to filtering as applied to the Kalman filter (KF). In what follows it will be understood that a comparison between the DLF approach and the Kalman filter is to be understood as the DLF approach applied to the Kalman filter and the classical Kalman filter.

Referring to (8) the posterior at any time $t_{n_m} \in T_O$ can be split up as follows

$$\pi(\boldsymbol{V}_{0:n_m}|\boldsymbol{Y}_{1:m}) \propto \pi(\boldsymbol{Y}_m|\boldsymbol{V}_{n_m})\pi(\boldsymbol{V}_{0:n_m}|\boldsymbol{Y}_{1:m-1}) \tag{9}$$

$$\pi(\boldsymbol{V}_{0:n_m}|\boldsymbol{Y}_{1:m-1}) = \pi(\boldsymbol{V}_{n_m}|\boldsymbol{V}_{n_m-1})\pi(\boldsymbol{V}_{0:n_m-1}|\boldsymbol{Y}_{1:m-1}). \tag{10}$$

Let $\boldsymbol{V}_{n|n}$, denote the posterior of $\boldsymbol{V}_n$ conditioned on observations up to $t_n$. If $t_n = t_{n_m} \in T_O$ for some $m$, then this distribution is simply $\pi(\boldsymbol{V}_{n_m}|\boldsymbol{Y}_{1:m})$, which by (9) can be written as

$$\begin{aligned}\boldsymbol{V}_{n_m|n_m} &\sim \int \pi(\boldsymbol{V}_{0:n_m}|\boldsymbol{Y}_{1:m})d\boldsymbol{V}_{n_m-1}\cdots d\boldsymbol{V}_0 \\ &\propto \pi(\boldsymbol{Y}_m|\boldsymbol{V}_{n_m})\int \pi(\boldsymbol{V}_{0:n_m-1}|\boldsymbol{Y}_{1:m-1})d\boldsymbol{V}_{n_m-1}\cdots d\boldsymbol{V}_0 \\ &= \pi(\boldsymbol{Y}_m|\boldsymbol{V}_{n_m})\pi(\boldsymbol{V}_{n_m}|\boldsymbol{V}_{n_m-1|n_m-1}).\end{aligned} \tag{11}$$

This prior $\pi(\boldsymbol{V}_{n_m}|\boldsymbol{V}_{n_m-1|n_m-1})$ can be calculated by applying the forward SDE solver (5) to $\boldsymbol{V}_{n-1|n-1}$, while the Likelihood $\pi(\boldsymbol{Y}_m|\boldsymbol{V}_{n_m})$ is determined entirely through the measurement errors. If, on the other hand, $t_n \notin T_O$, we simply have $\boldsymbol{V}_{n|n} \sim \pi(\boldsymbol{V}_n|\boldsymbol{V}_{n-1|n-1})$, which can be calculated through the model. Thus the posterior of $\boldsymbol{V}_{0:N}$ can be calculated sequentially for one $\boldsymbol{V}_n$ at a time, based on the posteriors up to the respective previous time step. All priors and likelihoods here, and therefore the posteriors too, are normally distributed.

We denote by $\langle \boldsymbol{V}_{n|n} \rangle$ and $\boldsymbol{P}_{n|n}$ the mean and covariance of $\boldsymbol{V}_{n|n}$ and by $\langle \boldsymbol{V}_{n|n-1} \rangle$ and $\boldsymbol{P}_{n|n-1}$ the mean and covariance of $\boldsymbol{V}_{n|n-1} \sim \pi(\boldsymbol{V}_n|\boldsymbol{V}_{n-1|n-1})$. Further, let $\boldsymbol{R}_m = \langle \boldsymbol{\epsilon}_m \boldsymbol{\epsilon}_m^\top \rangle$.

The KF produces sequential estimates for $\langle \boldsymbol{V}_{n|n} \rangle$ and $\boldsymbol{P}_{n|n}$, and thus for the posterior distribution, in two steps. In the *forecast* step the model is used to produce an initial estimate of $\pi(\boldsymbol{V}_n|\boldsymbol{V}_{n-1|n-1})$. Since the model is linear and the noise is (unbiased) normal, the prior at $t_n$ is estimated through the model and the posterior at the previous time step:

$$\langle \boldsymbol{V}_{n|n-1} \rangle = \boldsymbol{L}_{n-1} \langle \boldsymbol{V}_{n-1|n-1} \rangle + \Delta t \boldsymbol{f}_{n-1}, \quad n = 1, \ldots, N, \tag{12}$$

$$\boldsymbol{P}_{n|n-1} = \boldsymbol{L}_{n-1}\boldsymbol{P}_{n-1|n-1}\boldsymbol{L}_{n-1}^\top + \boldsymbol{Q}_{n-1}, \quad n = 1, \ldots, N, \tag{13}$$

with $\boldsymbol{Q}_{n-1} = \langle \boldsymbol{w}_{n-1}\boldsymbol{w}_{n-1}^\top \rangle$. Initial data is assumed to be known or a sample of a known distribution. If no observations are available at time $t_n$, the posterior is not affected by the likelihood, and thus $\langle \boldsymbol{V}_{n|n} \rangle = \langle \boldsymbol{V}_{n|n-1} \rangle$, and $\boldsymbol{P}_{n|n} = \boldsymbol{P}_{n|n-1}$. If, on the other hand, observations are available, i.e. $n = n_m$ with $t_n = t_{n_m} \in T_O$, an *analysis* step is performed, which takes in the mean and covariance of the prior $\langle \boldsymbol{V}_{n|n-1} \rangle$ and $\boldsymbol{P}_{n|n-1}$ and the mean and covariance of the likelihood $\langle \boldsymbol{Y}_n \rangle$ and $\boldsymbol{R}_n$ and calculates the moments of the posterior. For any step $t_{n_m}$ with observations, thus, the analysis step consists of the update

$$\langle \boldsymbol{V}_{n_m|n_m} \rangle = \langle \boldsymbol{V}_{n_m|n_m-1} \rangle + \boldsymbol{K}_{n_m} \left( \langle \boldsymbol{Y}_m \rangle - \boldsymbol{H}(\boldsymbol{y}_m) \langle \boldsymbol{V}_{n_m|n_m-1} \rangle \right), \tag{14}$$

$$\boldsymbol{P}_{n_m|n_m} = (\boldsymbol{I} - \boldsymbol{K}_{n_m}\boldsymbol{H}(\boldsymbol{y}_m))\boldsymbol{P}_{n_m|n_m-1}. \tag{15}$$

Here $\boldsymbol{H}(\cdot)$ evaluated at the vector $\boldsymbol{v} \in \Omega^I$ is the interpolation matrix defined as $\boldsymbol{H}(\boldsymbol{v}) := \left(\boldsymbol{H}(v^i)\right)_{i=1}^I \in \mathbb{R}^{I \times K}$. The term $\left(\langle \boldsymbol{Y}_m \rangle - \boldsymbol{H}(\boldsymbol{y}_m)\langle \boldsymbol{V}_{n_m|n_m-1} \rangle\right)$ in (14) is called the *innovation*. In (15) $\boldsymbol{I}$ is the $N$-dimensional identity matrix; $\boldsymbol{K}_{n_m}$ is called the Kalman gain and is defined as

$$\begin{aligned}\boldsymbol{K}_{n_m} =& \boldsymbol{P}_{n_m|n_m-1}\boldsymbol{H}(\boldsymbol{y}_m)^\top \\ &\cdot \left[\boldsymbol{H}(\boldsymbol{y}_m)\boldsymbol{P}_{n_m|n_m-1}\boldsymbol{H}(\boldsymbol{y}_m)^\top + \boldsymbol{R}_m\right]^{-1}.\end{aligned} \tag{16}$$

## 3. The Dynamic Likelihood Approach

In hyperbolic dynamics, (waves) information (along with uncertainties) will flow along characteristics. The DLF approach exploits the wave dynamics to propose a richer likelihood than other traditional problems. An observation $(\boldsymbol{y}_m, \boldsymbol{Y}_m)$ measured at time $t_{n_m} \in T_O$ is used to generate *pseudo-observations* $(\mathcal{F}_n\boldsymbol{y}_m, \mathcal{F}_n\boldsymbol{Y}_m)$, at times $t_n \in T$ with $t_n > t_{n_m}$. These are used to generate likelihood distributions in between observations. In fact, it is possible to produce likelihoods in the future, so that in principle, it is possible to do Bayesian estimation *in the future*. We refer to a likelihood, constructed from observations as well as from pseudo-observations, as the *dynamic likelihood*. The pseudo-observations are tightly coupled to the inherent model dynamics; in Foster and Restrepo (2022) we show how these are constructed for hyperbolic problems. Here we show how these could be formulated for advection-diffusion dynamics.

### a. Pseudo-Observations

The pseudo-observation at time $t_n$, from the observation $\boldsymbol{y}_m$ at time $t_{n_m}$, have locations $\mathcal{F}_n\boldsymbol{y}_m = \boldsymbol{x}_m(t_n) \in [0, L]^I$ which solve

$$
\begin{aligned}
dx_m^i(t) &= -c(x_m^i(t), t)dt, \qquad && t_{n_m} < t \le t_N, \\
x_m^i(t_{n_m}) &= y_m^i && i = 1, \dots, I
\end{aligned} \tag{17}
$$

A schematic of how $\mathcal{F}_n\boldsymbol{y}_m$ arises from $\boldsymbol{y}_m$ is depicted in Figure 1. The value $\mathcal{F}_n\boldsymbol{Y}_m$ of the pseudo-observations approximate $u$ at the location of the characteristics, and they include a measurement error. Namely, $\mathcal{F}_n\boldsymbol{Y}_m = \left(v(x_m^i(t_n), t_n) + \zeta^i(t_n)\right)_{i=1}^I$, where $\boldsymbol{\zeta}^i(t_m)$ has the same distribution as $\epsilon_m^i$.

Ideally, $\mathcal{F}_n\boldsymbol{Y}_m$ would be the values at $t = t_n$ of solutions $u^i(t) := u(x_m^i(t), t)$ for $t_{n_m} < t \le t_N$ of the following system of equations derived from (1) and (17),

$$
\begin{aligned}
du^i(t) =& \left( (D + \frac{1}{2}A^2)u_{xx}^i(t) + f(x_m^i(t), t) \right) dt \\
&+ BdW^u + Au_x^i(t)dW^c \\
u^i(t_{n_m}) =& Y_m^i.
\end{aligned} \tag{18}
$$

A schematic of how a single pseudo-observation $\mathcal{F}_n\boldsymbol{Y}_m^i$ would behave along a characteristic is depicted in Figure 2.

As the truth $u$ is assumed unknown, we recast equation (18) with an approximate form where the terms $u_x^i$ and $u_{xx}^i$ are replaced by approximations $v_{(x)}(x,t) \approx u_x(x,t)$ and $v_{(xx)}(x,t) \approx u_{xx}(x,t)$ which we refer to as first and second derivatives of the model. (We could obtain $v_{(x)}(x, t_n)$ and $v_{(xx)}(x, t_n)$ from the interpolated model $v$, this is not a requirement and they could just as well be given from another model or data). With these replacements (18)

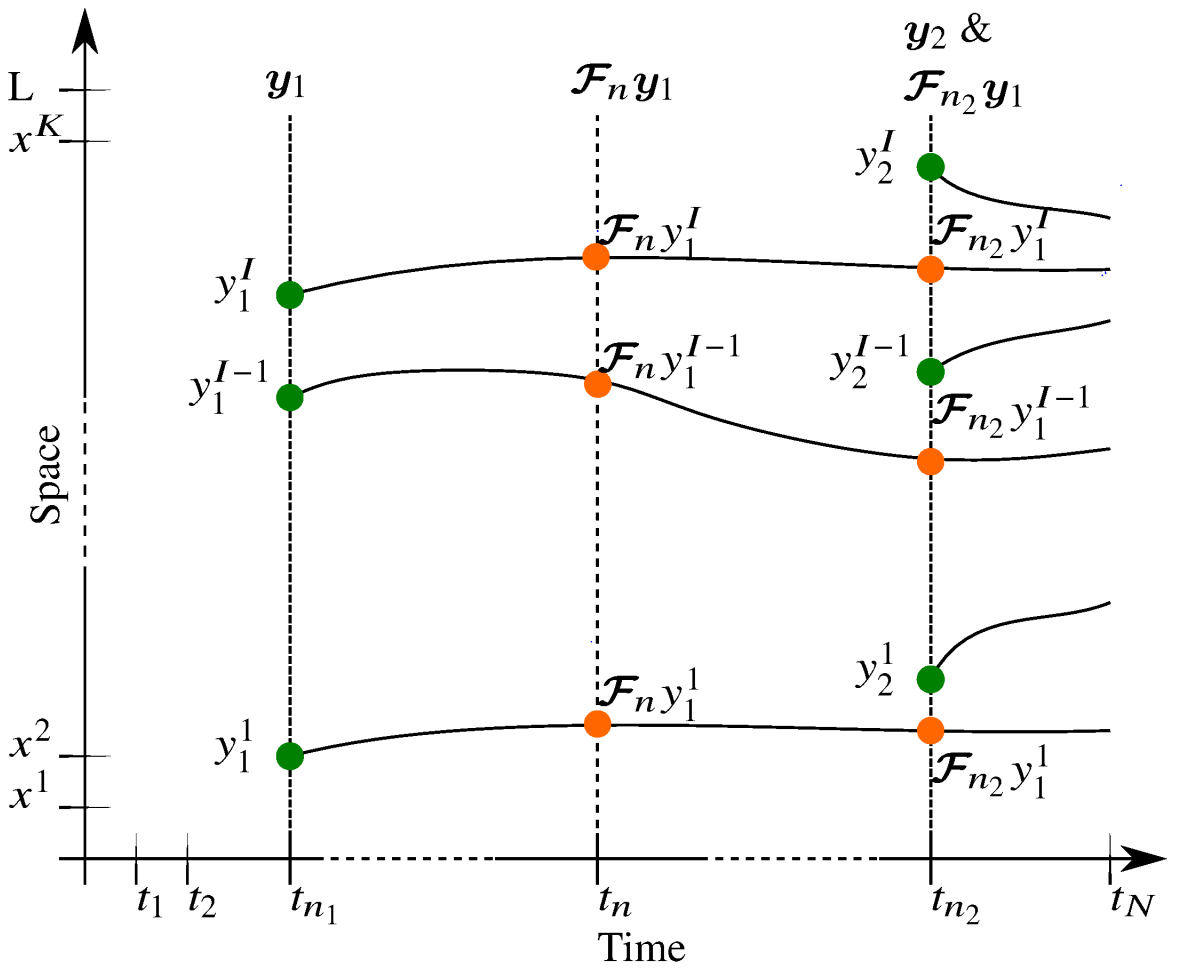

FIG. 1: Schematic depiction of how locations of pseudo-observations $\mathcal{F}_n\boldsymbol{y}_m$(orange) are derived from the locations of observations $\boldsymbol{y}_m$(green) by propagating along characteristics(black). In this graphic at $t_{n_1}$, only observations are available, thus only a classical filtering step can be performed, at $t_n$ only pseudo-observations are available, thus a regular DLF step can be performed, and at $t_{n_M}$ observations and pseudo-observations are available, which means an MDLF step would be performed. (See section a for the definition of DLF/MDLF-step and calculations of $\mathcal{F}_n\boldsymbol{y}_m$ as well as $\mathcal{F}_n\boldsymbol{Y}_m$.)

becomes

$$
\begin{aligned}
d\tilde{u}^i(t) =& \left( \left( D + \frac{A^2}{2} \right) v_{(xx)}(x_m^i(t), t) + f \right) dt \\
&+ Av_{(x)}((x_m^i(t), t))dW^c + BdW^u \\
\tilde{u}^i(t_{n_m}) =& Y_m^i,
\end{aligned} \tag{19}
$$

Provided $v_{(x)}$ and $v_{(xx)}$ along the characteristic $\boldsymbol{x}_m(t)$, equation (19) is a well-posed approximation of (17). We thus finally define the pseudo-observations

$$
\mathcal{F}_n\boldsymbol{y}_m := \left(x_m^i(t_n)\right)_{i=1}^I \tag{20}
$$

and

$$
\mathcal{F}_n\boldsymbol{Y}_m := \left(\tilde{u}^i(t_n)\right)_{i=1}^I, \tag{21}
$$

where the $x_m^i(t)$ and $\tilde{u}^i(t)$ solve equations (17) and (19) respectively for $i = 1, \dots, I$. At estimation times $t_n \in T$, the pseudo-observations $(\mathcal{F}_n\boldsymbol{y}_m, \mathcal{F}_n\boldsymbol{Y}_m)$ are treated like real observations, and the *pseudo-observation error* $\boldsymbol{\zeta}_n \in \mathbb{R}^I$ is defined as

$$
\boldsymbol{\zeta}_n = \mathcal{F}_n\boldsymbol{Y}_m - v(\mathcal{F}_n\boldsymbol{y}_m, t_n). \tag{22}
$$

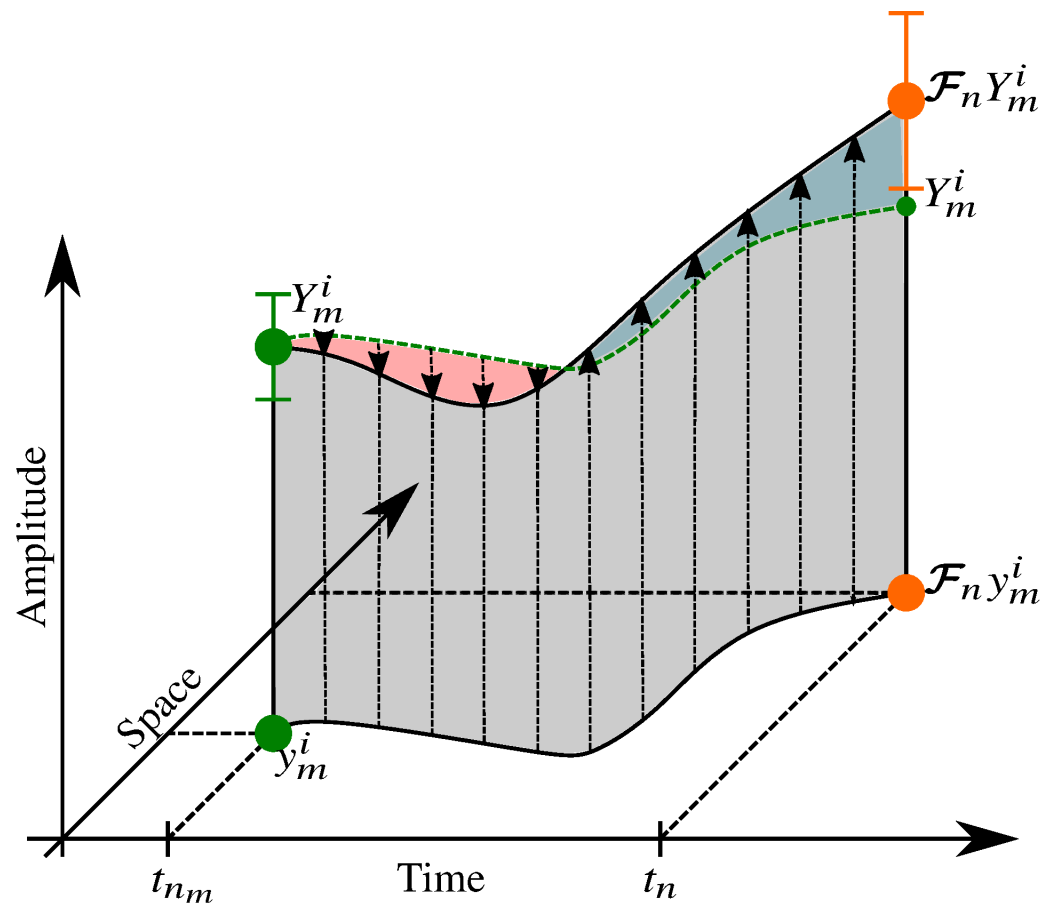

FIG. 2: Schematic depiction of how values of pseudo-observations $\mathcal{F}_n \boldsymbol{Y}_m$ are derived by evolving observations $\boldsymbol{Y}_m$ along characteristics. At $t_{n_m}$, we see a single observation $Y_m^i$ at location $y_m^i$(green dots). At $t_n$ the pseudo-observation $\mathcal{F}_n Y_m^i$ at location $\mathcal{F}_n y_m^i$ is depicted(orange dots.) As a dashed green line, we see a curve of constant amplitude $Y_m^i$ along the characteristic starting at $y_m^i$. From $Y_m^i$ to $\mathcal{F}_n Y_m^i$ in black the actual trajectory of the pseudo-observation is shown. Underlayed in red are times, when the pseudo-observation is smaller than $Y_m^i$, and underlayed in blue are times, at which the pseudo-observation exceeds the value of $Y_m^i$.

$\zeta_n$ is assumed to be distributed according to a mean zero Gaussian distribution. Recall that the value of $v(\mathcal{F}_n \boldsymbol{y}_m, t_n)$ can be obtained from $\boldsymbol{V}_n$ using (6). Equation (22) thus induces the distribution $\pi(\mathcal{F}_n \boldsymbol{Y}_m | \boldsymbol{V}_n)$, which will be used to inform the Likelihood at times $t_n$ used in the DLF and discussed in section b. For later reference, we define $\mathcal{F}_n \boldsymbol{y}_{\ell:m}$ and $\mathcal{F}_n \boldsymbol{Y}_{\ell:m}$ for $n_\ell < n_m < n$ as, respectively, the column vectors obtained by concatenating all $\mathcal{F}_n \boldsymbol{y}_{n_{m'}}$ and $\mathcal{F}_n \boldsymbol{Y}_{n_{m'}}$ with $t_{n_{m'}} \in T_O$ and $t_{n_\ell} \le t_{n_{m'}} \le t_{n_m}$. Namely,

$$\mathcal{F}_n \boldsymbol{y}_{\ell:m} = \begin{pmatrix} \mathcal{F}_n \boldsymbol{y}_\ell \\ \vdots \\ \mathcal{F}_n \boldsymbol{y}_m \end{pmatrix}, \qquad \mathcal{F}_n \boldsymbol{Y}_{\ell:m} = \begin{pmatrix} \mathcal{F}_n \boldsymbol{Y}_\ell \\ \vdots \\ \mathcal{F}_n \boldsymbol{Y}_m \end{pmatrix}$$

Thus $(\mathcal{F}_n \boldsymbol{y}_{\ell:m}, \mathcal{F}_n \boldsymbol{Y}_{\ell:m})$ contains all pseudo-observations at time $t_n$ derived from observations that were measured at times between $t_{n_\ell}$ and $t_{n_m}$.

### *b. Formulation of the Filter*

We define for the duration of this section $n \in \{1, \ldots, N\}$ and $m \in \{1, \ldots, M\}$ such that $n_m \le n \le n_{m+1}$. Under the assumption posed on the model, observations and pseudo-observations, the DLF yields an alternative form of the posterior of the model, up to time $t_n$, conditioned on observations measured up to time $t_{n_m}$:

$$\pi(\boldsymbol{V}_{0:n} | \boldsymbol{Y}_{1:m}) \propto \pi(\boldsymbol{Y}_m, \mathcal{F}_n \boldsymbol{Y}_{1:m-1} | \boldsymbol{V}_n) \pi(\boldsymbol{V}_{0:n} | \boldsymbol{Y}_{1:m-1}). \tag{23}$$

Let $\boldsymbol{V}_{n|n}$ denote the posterior of $\boldsymbol{V}_n$ conditioned on observations up to time $t_n$. It is distributed as

$$\begin{aligned} \boldsymbol{V}_{n|n} &\sim \pi(\boldsymbol{V}_n | \boldsymbol{Y}_m, \mathcal{F}_n \boldsymbol{Y}_{1:m-1}) \\ &= \int \pi(\boldsymbol{V}_{0:n_m} | \boldsymbol{Y}_{1:m}) d\boldsymbol{V}_{n_m-1} \cdots d\boldsymbol{V}_0 \\ &\propto \pi(\boldsymbol{Y}_m, \mathcal{F}_n \boldsymbol{Y}_{1:m-1} | \boldsymbol{V}_{n_m}) \\ &\quad \cdot \int \pi(\boldsymbol{V}_{0:n_m-1} | \boldsymbol{Y}_{1:m-1}) d\boldsymbol{V}_{n_m-1} \cdots d\boldsymbol{V}_0 \\ &= \pi(\boldsymbol{Y}_m, \mathcal{F}_n \boldsymbol{Y}_{1:m-1} | \boldsymbol{V}_n) \pi(\boldsymbol{V}_n | \boldsymbol{V}_{n-1|n-1}). \end{aligned} \tag{24}$$

We therefore see that this posterior at time $t_n$ is entirely determined by the likelihood of observations and pseudo-observations at $t_n$ conditioned on the model up to this point $\pi(\boldsymbol{Y}_m, \mathcal{F}_n \boldsymbol{Y}_{1:m-1} | \boldsymbol{V}_n)$, and the distribution of the model conditioned on the posterior at the previous time step $\pi(\boldsymbol{V}_n | \boldsymbol{V}_{n-1|n-1})$. The latter is used as a prior at $t_n$.

This shows that, as in the KF approach and given the observations and pseudo-observations, the distribution of $\boldsymbol{V}_{n|n}$ can be found sequentially based on the posterior distribution at the previous time step $\boldsymbol{V}_{n-1|n-1}$. For $t_n \le \min(T_O)$ there are no pseudo-observations, i.e. $\mathcal{F}_n \boldsymbol{Y}_{1:n-1} = (\ )$ is an empty vector, thus the KF approach is recovered.

If on the other hand $t_n > \min(T_O)$ there are two cases to be considered:

(i) there are no observations at $t_n \notin T_O$ and thus $\pi(\boldsymbol{Y}_m, \mathcal{F}_n \boldsymbol{Y}_{1:m-1} | \boldsymbol{V}_n) = \pi(\mathcal{F}_n \boldsymbol{Y}_{1:m-1} | \boldsymbol{V}_n)$,

(ii) $t_n = t_{n_m} \in T_O$ and therefore $\pi(\boldsymbol{Y}_m, \mathcal{F}_{n_m} \boldsymbol{Y}_{1:m-1} | \boldsymbol{V}_{n_m}) = \pi(\mathcal{F}_{n_m} \boldsymbol{Y}_{1:m-1} | \boldsymbol{V}_{n_m}) \pi(\boldsymbol{Y}_m | \boldsymbol{V}_{n_m})$.

We call a step of the DLF algorithm in the first case a *DLF update* and in the second case a *multi analysis DLF (MDLF) update*. In the first case one proceeds as follows:

The DLF update:

$$\begin{aligned} \text{Predict:} \langle \boldsymbol{V}_{n|n-1} \rangle &= \boldsymbol{L}_{n-1} \langle \boldsymbol{V}_{n-1|n-1} \rangle + \Delta t f_{n-1} \\ \boldsymbol{P}_{n|n-1} &= \boldsymbol{L}_{n-1} \boldsymbol{P}_{n-1|n-1} \boldsymbol{L}_{n-1}^\top + \boldsymbol{Q}_{n-1} \\ (\boldsymbol{y}_\mathcal{F}, \boldsymbol{Y}_\mathcal{F}) &= (\mathcal{F}_n \boldsymbol{y}_{1:m}, \mathcal{F}_n \boldsymbol{Y}_{1:m}) \\ \boldsymbol{R}_\mathcal{F} &= \mathrm{Cov}(\boldsymbol{Y}_\mathcal{F}, \boldsymbol{Y}_\mathcal{F}) \\ \text{Analysis:} \langle \boldsymbol{V}_{n|n} \rangle &= \boldsymbol{V}_{n|n-1} + \boldsymbol{K}_n \left( \langle \boldsymbol{Y}_\mathcal{F} \rangle - \boldsymbol{H}(\boldsymbol{y}) \langle \boldsymbol{V}_{n|n-1} \rangle \right) \\ \boldsymbol{P}_{n|n} &= (\boldsymbol{I} - \boldsymbol{K}_n \boldsymbol{H}(\boldsymbol{y}_\mathcal{F})) \boldsymbol{P}_{n|n-1} \end{aligned}$$

where

$$\boldsymbol{K}_n = \boldsymbol{P}_{n|n-1} \boldsymbol{H}(\boldsymbol{y}_\mathcal{F})^\top \left( \boldsymbol{H}(\boldsymbol{y}_\mathcal{F}) \boldsymbol{P}_{n|n-1} \boldsymbol{H}(\boldsymbol{y}_\mathcal{F})^\top + \boldsymbol{R} \right)^{-1}$$

This covers the case, where only pseudo-observations are available at time $t_n$. Derivations of this mirror the derivation of the KF exactly, if observations are replaced by pseudo-observations. If both real and pseudo-observations are available at $t_n = t_{n_m}$, i.e. $n = n_m$, a multi-analysis step is performed. The prediction step is identical as in the previous case. The Analysis step however changes as follows:

The MDLF update:

$$
\begin{aligned}
\text{Predict:}\langle \boldsymbol{V}_{n_m|n_m-1}\rangle &= \boldsymbol{L}_{n_m-1}\langle \boldsymbol{V}_{n_m-1|n_m-1}\rangle + \Delta t f_{n_m-1}\\
\boldsymbol{P}_{n_m|m_m-1} &= \boldsymbol{L}_{n_m-1}\boldsymbol{P}_{n_m-1|n_m-1}\boldsymbol{L}^\top_{n_m-1} + \boldsymbol{Q}_{n_m-1}\\
(\boldsymbol{y}_{\mathcal{F}}, \boldsymbol{Y}_{\mathcal{F}}) &= \left(\mathcal{F}_{n_m}\boldsymbol{y}_{1:m-1}, \mathcal{F}_{n_m}\boldsymbol{Y}_{1:m-1}\right)\\
\boldsymbol{R}_{\mathcal{F}} &= \operatorname{Cov}(\boldsymbol{Y}_{\mathcal{F}}, \boldsymbol{Y}_{\mathcal{F}})\\
\text{Update:}\langle \boldsymbol{V}_{n_m|n_m}\rangle &= \langle \boldsymbol{V}_{n_m|n_m-1}\rangle\\
&\quad + \boldsymbol{K}^*_{n_m}(\langle \boldsymbol{Y}_m\rangle - \boldsymbol{H}(\boldsymbol{y}_m)\langle \boldsymbol{V}_{n|n-1}\rangle)\\
&\quad + \boldsymbol{J}_{n_m}\left(\langle \boldsymbol{Y}_{\mathcal{F}}\rangle - \boldsymbol{H}(\boldsymbol{y}_{\mathcal{F}})\langle \boldsymbol{V}_{n|n-1}\rangle\right)\\
\boldsymbol{P}_{n_m|n_m} &= (\boldsymbol{I} - \boldsymbol{K}^*_{n_m}\boldsymbol{H}(\boldsymbol{y}_m)\\
&\quad - \boldsymbol{J}_{n_m}\boldsymbol{H}(\boldsymbol{y}_{\mathcal{F}}))\boldsymbol{P}_{n_m|n_m-1}
\end{aligned}
$$

where

$$
\begin{aligned}
\boldsymbol{K}^*_{n_m} &= (\boldsymbol{I} - (\boldsymbol{I} - \boldsymbol{K}\boldsymbol{H}(\boldsymbol{y}_m))\boldsymbol{D}\boldsymbol{H}(\boldsymbol{y}_{\mathcal{F}}))\boldsymbol{K}\\
\boldsymbol{J}_{n_m} &= (\boldsymbol{I} - \boldsymbol{K}_{n_m}\boldsymbol{H}(\boldsymbol{y}_m))\boldsymbol{D}
\end{aligned}
$$

with

$$
\begin{aligned}
\boldsymbol{K} =& \boldsymbol{P}_{n_m|n_m-1}\boldsymbol{H}(\boldsymbol{y}_{\mathcal{F}})^\top\\
&\cdot(\boldsymbol{R}_{n_m} + \boldsymbol{H}(\boldsymbol{y}_m)\boldsymbol{P}_{n_m|n_m-1}\boldsymbol{H}(\boldsymbol{y}_m)^\top)^{-1}\\
\boldsymbol{D} =& \boldsymbol{P}_{n_m|n_m-1}\boldsymbol{H}(\boldsymbol{y}_{\mathcal{F}})^\top\left(\boldsymbol{R}_{\mathcal{F}} + \boldsymbol{H}(\boldsymbol{y}_{\mathcal{F}})\boldsymbol{P}_{n_m|n_m-1}\boldsymbol{H}(\boldsymbol{y}_{\mathcal{F}})^\top\right.\\
&\left.-\boldsymbol{H}(\boldsymbol{y}_{\mathcal{F}})\boldsymbol{K}_{n_m}\boldsymbol{H}(\boldsymbol{y}_m)\boldsymbol{P}_{n_m|n_m-1}\boldsymbol{H}(\boldsymbol{y}_{\mathcal{F}})^\top\right)^{-1}.
\end{aligned}
$$

The derivation of these gains can be found in either Restrepo (2013) or Foster and Restrepo (2022). Note that the size of the matrices to invert in the calculation of these gains is determined by the size of $\boldsymbol{R}_{\mathcal{F}} \in \mathbb{R}^{mI\times mI}$. In practice, it might be reasonable to limit the size of this matrix by not using all pseudo-observations indefinitely, but rather discarding some in a trade-off between accuracy and complexity. As an example, replacing $\mathcal{F}_n\boldsymbol{y}_{1:m}$ and $\mathcal{F}_n\boldsymbol{Y}_{1:m}$ by $\mathcal{F}_n\boldsymbol{y}_{k_n:m}$ and $\mathcal{F}_n\boldsymbol{Y}_{k_n:m}$ respectively in the previous algorithms, for some $k_n$, would discard pseudo-observations derived from the oldest observations as time goes on, thus limiting the size of $\boldsymbol{R}_{\mathcal{F}}$ to $(m-k_n)I(m-k_n)I$. The overall complexity of the DLF and the effects of discarding "older" pseudo-observations like this are discussed in the next section.

## 4. Analysis of Computational Complexity of the DLF

The DLF approach has an added computational overhead, as compared to its standard counterpart. In what follows we estimate the overhead of the DLF approach to Kalman filtering, as compared to the native Kalman filter. As will be shown subsequently, improvements in the estimates obtained using DLF may offset the added computational burden. The additional computational cost of the DLF, compared to the model on its own, stems from two sources:

1. Calculating the pseudo-observations $(\mathcal{F}_n\boldsymbol{y}_{1:m-1}, \mathcal{F}_n\boldsymbol{Y}_{1:m-1})$.
2. Calculating and applying the gains $\boldsymbol{K}_n, \boldsymbol{K}^*_n, \boldsymbol{J}_n$.

If an explicit ODE solver is used to solve equations (17) and (19) to calculate $(\mathcal{F}_n\boldsymbol{y}_m, \mathcal{F}_n\boldsymbol{Y}_m)$ at $t_{n_m} \in T_O$ the complexity is linear in the number of time steps $\frac{t_n - t_{n_m}}{\Delta t}$. Since $\mathcal{F}_n\boldsymbol{y}_m \in \mathbb{R}^I$ and $\mathcal{F}_n\boldsymbol{Y}_m$ is normally distributed in $\mathbb{R}^I$ these calculations over a single time step are dominated by the calculation of the covariance $\boldsymbol{R}_{n,m} = \operatorname{Cov}(\mathcal{F}_n\boldsymbol{Y}_m, \mathcal{F}_n\boldsymbol{Y}_m) \in \mathbb{R}^{I\times I}$ leading to an overall complexity of order $O\left(\frac{t_n - t_{n_m}}{\Delta t} I^2\right)$ to calculate $(\mathcal{F}_n\boldsymbol{y}_m, \mathcal{F}_n\boldsymbol{Y}_m)$ for a single $t_{n_m} \in T_O$. Considering this has to be done for all pseudo-observations $(\mathcal{F}_N\boldsymbol{y}_{1:M-1}, \mathcal{F}_N\boldsymbol{Y}_{1:M-1})$ this leads to a complexity of $O\left(\sum_{m=1}^M \frac{t_N - t_{n_m}}{\Delta t} I^2\right) \leq O\left(MNI^2\right)$.

The complexity of calculating the gains at time $t_n \in T$ is dominated by inverting matrices of the same size as the covariance matrix $\boldsymbol{R}_{\mathcal{F}}$. At time $t_n$ let this size be $s_n \times s_n$. Note that for times $t_{n_m} \leq t_n < t_{n_{m+1}}$ this size is $s_n = mI$. The computational cost of matrix inversion is cubic in the number of rows and therefore the complexity for the $n-$th time step is $O(I^3M^3)$. Over all time steps this is bounded by $O(NM^3I^3)$. The complexity of the KF arises from inverting matrices of size $I \times I$ at each of the $M$ observation times, thus having an overall complexity of $O(MI^3)$. The DLF's overall complexity is $O(MI^3 + MNI^2)$.

Note however, that in all these complexity estimates so far we assumed that all pseudo-observations $(\mathcal{F}_n\boldsymbol{y}_m, \mathcal{F}_n\boldsymbol{Y}_m)$ are used for all $t_n > t_{n_m}$. As $t_n - t_{n_m}$ increases so does the uncertainty associated with $(\mathcal{F}_n\boldsymbol{y}_m, \mathcal{F}_n\boldsymbol{Y}_m)$. Therefore discarding this pseudo-observation eventually would have only a small effect on overall accuracy. Thus if run time is of the essence, discarding some pseudo-observations after they outlived their usefulness can help to keep complexity in check. One way to proceed with this program is to set an upper bound for the number of pseudo-observations assimilated at any given time, or by setting a threshold on uncertainty, discarding the relatively most uncertain observations/pseudo-observations.

As an example, let us assume the number of pseudo-observations is capped at an integer multiple of $I$, say $pI$, and the oldest pseudo-observations are discarded every time this threshold is reached. This means at any given time $t_n$ only observations from the last $p$ observation times are used to calculate pseudo-observations. At an observation $t_n = t_{n_m}$, a multi-analysis step is performed, before the oldest pseudo-observations $(\mathcal{F}_{n_m}\boldsymbol{y}_{m-p}, \mathcal{F}_{n_m}\boldsymbol{Y}_{m-p})$ from time $t_{n_{m-p}}$ are discarded. In other words $(\mathcal{F}_n\boldsymbol{y}_{1:m-1}, \mathcal{F}_n\boldsymbol{Y}_{1:m-1})$ in the DLF algorithm is replaced by $(\mathcal{F}_n\boldsymbol{y}_{m-p:m-1}, \mathcal{F}_n\boldsymbol{Y}_{m-p:m-1})$. This limits the number of pseudo-observations concurrently in the algorithm to $pI$ from the previous maximum of $MI$. We get new complexity estimates under these circumstances by making the replacement $M \longrightarrow p$ in the previous estimates. Thus the complexity of the DLF falls to $O(Np^3I^3 + pNI^2)$. Assuming $p$ is picked reasonably small ($p \ll N$ and $p \ll I$) this reduces further to $O(NI^3)$, which is comparable to a standard KF's complexity of be $O(MI^3)$. (While the calculation of the gains can be a bottleneck of the algorithm, it is noteworthy that as long as neither $c$ nor $f$ in equations (2) and (3) depend on the observations $\boldsymbol{Y}_{1:M}$, the gains are also independent of $\boldsymbol{Y}_{1:M}$ and can thus be calculated offline).

## 5. Numerical Results and Comparisons

We contrast the DLF approach, as applied to the KF (the DLF), with estimates obtained by the native KF (the KF). We will also discuss the numerical details of our implementation and introduce the metrics, based on which we will evaluate the performance of the two approaches.

### *a. Computational Details*

Since the DLF approach was developed specifically for hyperbolic (wave) equations, we are especially interested in determining to what extent it can handle advection as well as macroscale diffusion. We thus introduce the nondimensional quantity $\alpha := D/c_0L$, where $D$ is the diffusion coefficient, $c_0$ is the typical size of the velocity $C$, and $L$ is the characteristic length which is taken as the length of the domain, to capture the extent to which diffusion processes and advection qualitatively affect the solution. Let the primed nondimensional quantities be

$$x' = \frac{x}{L}, \quad t' = \frac{tc_0}{L}, \quad u' = \frac{u}{u_0}, \quad F' = \frac{FL}{u_0c_0}, \quad C' = \frac{C}{c_0},$$

then Eq. (1) is now recast as

$$\begin{aligned} u_t - (c+\chi)u_x &= \alpha u_{xx} + f + \phi, \quad t > 0, x \in [0,1], \\ u(x,0) &= \mathcal{U}(x), \quad x \in [0,1], \end{aligned} \tag{25}$$

in dimensionless units, having dropped the primes. It is still assumed that $\phi$ and $\chi$ are noise terms generated by Wiener processes. Therefore $\chi dt = AdW^c$ and $\phi dt = BdW^u$ still hold.

We run all numerical examples shown in this section on a domain $[0,1]$ with periodic boundary conditions and times spanning $t_0 = 0$ to $t_N = 0.5$. We chose $\Delta x = 0.01$ and $\Delta t = 0.005$, in dimensionless units for the discretization. The wave velocity is set to be

$$C(x,t)dt = \cos(5\pi t)\,dt + AdW^c(x,t) + \tilde{A}d\tilde{W}_c(t)$$

Note that we split the wave noise term into $AdW^c(x,t) + \tilde{A}d\tilde{W}_c(t)$. Both $dW^c$ and $d\tilde{W}_c$ are incremental Wiener processes, but $dW^c$ will be assumed to be uncorrelated in space, *i.e.*, $\langle dW^c(x)dW^c(y)\rangle = \delta_{x,y}$ for all $x, y \in X$, while $d\tilde{W}_c$ is independent of $x$. We will set $A = 0.05$ for all numerical experiments but will consider two cases of $\tilde{A}$ when simulating the truth. In the first case $\tilde{A} = 0$, while in the second $\tilde{A} = 1$. The model will be unaware of $\tilde{A}$, i.e. assume $\tilde{A} = 0$ in both cases. This introduces a systematic error in $C$ for the model that the KF and DLF will have to overcome. The effects of this will be discussed in section 3.

In all of the following examples, we assume that the forcing is given by

$$F(x,t)dt = BdW^u(x,t) = 0.05dW^u(x,t)$$

where $dW^u(x,t)$ is an incremental Wiener process in time for each $x$ and $\langle dW^u(x,t)dW^u(y,t)\rangle = \delta_{x,y}$ for all $x, y \in X$.

The initial data will be

$$u_0(x) = \sigma\exp(-250(x-\theta)^2).$$

The amplitude $\sigma$ and phase $\theta$ will be deterministic or chosen from random distributions. We will discuss three kinds of initial data for the model: (i) the deterministic case where $\sigma = 1$ and $\theta = \frac{1}{2}$ (see Section 1); (ii) with uncertain amplitude: $\sigma \sim \mathcal{U}\left[\frac{1}{2}, \frac{3}{2}\right]$, where $\mathcal{U}$ indicates a uniform distribution, and $\theta = \frac{1}{2}$; and (iii) the case of $\sigma = 1$ and uncertain phase $\theta \sim \mathcal{U}[0,1]$. In all these cases the amplitude and phase of the truth are fixed to $\sigma = 1$ and $\theta = \frac{1}{2}$. Cases (ii) and (iii) will be highlighted in Section 2.

The *truth* will be used to test the outcomes as well as to generate observations. The truth is computed through Strang-splitting Strang (1968). Equation (1) is split into a noisy advection equation

$$u_t - (c+\chi)u_x = f + \phi$$

and a deterministic diffusion equation

$$u_t = \alpha u_{xx},$$

which are then used to generate a solution sequentially. The diffusive step is solved via FFT by calculating an exact solution in Fourier space. For the advective part of the splitting, we chose a Lax-Wendroff scheme in space and a stochastic Runge-Kutta scheme in time. The chosen RK method is a third-order scheme with second-order weak convergence Rossler (2009).

The *model* is given to us as a first-order split step procedure, splitting (1) into the same noisy advection equation and diffusion equation as for the truth. The diffusion step is, again, solved exactly via FFT. The advective step of the model is solved with an upwind scheme in space and an explicit Euler scheme in time. The Courant number for all methods is 1.

*Observations* are derived from the truth. These are made available at observation times $t_{n_m} \in \{0.05, 0.1, 0.15, ..., 0.45\} = T_O$ by uniformly drawing their locations $y_m^1, \dots, y_m^I \in X$ from the grid, without repetition. Their corresponding values are then determined from the truth via $Y_m^i = u(y_m^i, t_{n_m}) + \epsilon_m^i$, where the $\epsilon_m^i$ are drawn independently from a mean zero normal distribution with variance $10^{-4}$. Throughout the next sections, the number of observations at each observation time will take values $I \in \{10, 20, 40, 60\}$.

Both the KF and the DLF require the interpolation operator $\boldsymbol{H}(x)$. In these tests, we use a simple linear interpolation operator, which for $x \in [0, L]$ is defined as

$$\boldsymbol{H}(x) := \big((1 - r(x)\delta_{k,\ell}) + r(x)\delta_{k,\ell+1}\big)_{k=0}^{K} \in \mathbb{R}^{1\times K}. \quad (26)$$

Here, $r(x) = \mod(x, \Delta x)$ is the remainder of $\frac{x}{\Delta x}$, $\delta_{i,j}$ is the Kronecker delta, and $\ell \in \{1, \dots, K\}$ is chosen such that $x^\ell \le x < x^{\ell+1}$ for grid points $x^\ell, x^{\ell+1} \in X$. (To account for periodic boundary conditions, assume $x^K < x^0 = L$ for $x^K \le x < L$).

To calculate $\mathcal{F}_n \boldsymbol{y}_m$, we use an explicit Euler algorithm to solve equation (17). Since for $\mathcal{F}_n \boldsymbol{Y}_m$, the calculation of the mean and covariance is sufficient for the algorithm we use an explicit Euler algorithm based on equation (19) for each of those as well. The derivatives of the model are therefore only required at $t_n \in T$. We obtain them through the model as $v_{(x)}(x, t_n) = \boldsymbol{H}(x)\nabla_x \boldsymbol{V}_n$ and $v_{(xx)}(x, t_n) = \boldsymbol{H}(x)\nabla_{xx}\boldsymbol{V}_n$, where $\nabla_x$ and $\nabla_{xx}$ are center difference operators approximating first and second derivatives. Thus over a single time step from $t_{n-1} \ge t_m$ to $t_n$ and with $\tilde{\boldsymbol{H}}_n := \boldsymbol{H}(\mathcal{F}_n \boldsymbol{y}_m)$ and $\boldsymbol{R}_{n,m} = \mathrm{Cov}(\mathcal{F}_n \boldsymbol{Y}_m, \mathcal{F}_n \boldsymbol{Y}_m)$ we use

$$\begin{aligned}\langle \mathcal{F}_n \boldsymbol{Y}_m \rangle &= \langle \mathcal{F}_{n-1} \boldsymbol{Y}_m \rangle \\ &\quad + \Delta t \cdot \left(D + \frac{A^2}{2}\right) \tilde{\boldsymbol{H}}_{n-1} \nabla_{xx} \langle \boldsymbol{V}_{n-1|n-1} \rangle, \end{aligned} \quad (27)$$

and

$$\begin{aligned}\boldsymbol{R}_{n,m} &= \boldsymbol{R}_{n-1,m} + \Delta t \cdot \\ &\left(B^2 \boldsymbol{I} + A^2 (\nabla_x \tilde{\boldsymbol{H}}_{n-1} \langle V_{n-1|n-1} \rangle)(\nabla_1 \tilde{\boldsymbol{H}}_{n-1} V_{n-1|n-1})^\top\right) \\ &+ \Delta t^2 \cdot \left(D + \frac{A^2}{2}\right) \tilde{\boldsymbol{H}}_{n-1} \nabla_{xx} \boldsymbol{P}_{n-1|n-1} \nabla_{xx}^\top \tilde{\boldsymbol{H}}_{n-1}^\top. \end{aligned} \quad (28)$$

The second line in (17) is an explicit Euler scheme modeling the system noise, while the third line tracks the errors introduced due to the uncertainty of $\boldsymbol{V}_n$.

To quantitatively compare the traditional KF to the DLF approach to the KF we calculate the following metrics: the Residual Mean Square (RMS) error, the Mass error, the Center of Mass (CoM) error, and the probabilistic Calibration. These are given by

$$\text{RMS error:} \sqrt{\Delta t \Delta x \sum_{n=1}^{N} \sum_{k=1}^{K} |U^i(t_n) - \langle V_{n|n}^k \rangle)|^2}$$

$$\text{Mass error:} \sqrt{\Delta t \Delta x^2 \sum_{n=1}^{N} \left| \sum_{i=1}^{N} |U^i(t_n)| - \sum_{i=1}^{N} |\langle V_{n|n}^k \rangle| \right|^2}$$

$$\text{CoM error:} \sqrt{\Delta t \sum_{n=1}^{N} \left| \frac{\sum_{k=1}^{K} |U^k(t_n) x^k|}{\sum_{k=1}^{K} |U^k(t_n)|} - \frac{\sum_{k=1}^{K} |\langle V_{n|n}^k \rangle x^k|}{\sum_{k=1}^{K} |\langle V_{n|n}^k \rangle|} \right|^2}$$

$$\begin{aligned}\text{Calibration:} \frac{\Delta t}{t_N N} \sum_{n=1}^{N} \sum_{k=1}^{K} \mathbb{1}\Big((U^k(t_n) \\ - \langle V_{n|n}^k \rangle) < 2\sqrt{\mathrm{Var}(V_{n|n}^k)}\Big),\end{aligned}$$

where $\mathbb{1}(\cdot)$ is the indicator function.

The RMS error tracks the sum of local errors between model and truth, while the Mass error determines how accurately the total mass in the system is captured. The CoM error remains small if the position of the center of mass is captured well by the model. This is important in advection-dominated problems and in our examples will be mostly determined by how well the position of the maximum is captured over time. The Calibration measures the percentage of times the truth is within two of the estimated standard deviations of the model. If both noise and model error are normally distributed and captured correctly by the uncertainty, this value should be approximately 95%. Higher values indicate the variance is overestimated, while smaller values mean the uncertainty is larger than estimated. In the following sections, we will see these measures evaluated in total and at specific times. When evaluated at a specific time $t_n$, the $\Delta t$ and the sum over $n$ will be dropped and what remains will be evaluated at the corresponding $n$. To account for the randomness in the generation of the truth,(and the initial data of the model, where applicable,) all records of these four metrics from hereon will be based on 50 runs each: Line plots of RMS, Mass, Com error and Calibration will show their mean value over 50 runs, while all box plots will be based on their respective minimum, maximum and the 25%, 50% and 75% quantiles. Between each of these runs the truth, observations, and initial data will be regenerated. For comparability, the KF and DLF will use the same Observations and initial data for each individual run.

### *b. Comparing the KF, and the DLF Outcomes*

We will compare outcomes for deterministic and aleatoric initial data, uncertainties in the phase speed and as a function of $\alpha$.

#### 1) Comparing the KF and DLF in Systems with Deterministic Initial Data

We will compare posterior predictions $\boldsymbol{V}_{n|n}$. The comparison is conducted in a setting where both the noise of the phase speed $dW^c$ and the noise of the forcing $dW^u$ are uncorrelated in space, i.e., $\langle dW(x)dW(y)\rangle = \delta_{x,y}$. Both methods are provided with deterministic initial data, $\sigma = 1$, and $\theta = \frac{1}{2}$ and $\tilde{A} = 0$ . Throughout this comparison, both the relative diffusion $\alpha$ and the number of data points per observation $I$ will be varied. We will demonstrate that the DLF outperforms the KF in terms of RMS, Mass error, and Calibration, particularly when the number of data points is sparse and when $\alpha$ is small. The amplitude of the wave noise will be fixed at $A = 0.05$ for the entirety of this section.

We first examine an individual run in the advection-dominated case. For $\alpha = 0.01$, Figure 3 shows the truth (right), the model prediction through the KF (left), and the prediction of the DLF (middle). Both filters were presented with $I = 20$ data points per observation time. The location of these data points is randomly selected at each observation time but is identical between the KF and the DLF. Note that the trajectories of pseudo-observations depicted extend beyond the availability of observation, providing Bayesian predictions at times $t > 0.45$, which can be considered the future.

At the observation sites, both the KF and the DLF pick up the values of the observed data and adjust their predictions during the analysis step. The DLF manages to maintain these adjusted values over the simulation time, while in the KF, adjustments due to observations quickly vanish due to diffusion. This is not surprising, as the DLF reinforces the information gathered from observations through pseudo-observations that move along characteristics.

Figure 4 confirms that the DLF captures the truth more accurately, as quantified by the metrics. We examine the time series of these four metrics for the same parameters and conditions used to generate the previous example. Starting from time $t = 0.05$, the first observation time, there is a clear divide between the KF and the DLF in the RMS and Mass errors, with the DLF performing significantly better. The same can be said for the Calibration, though the advantage of the DLF is less pronounced. There appears to be no such clear trend for the CoM error.

As the last part of this set of experiments, we examine if these trends hold up when the remaining parameters $\alpha$ and $I$ are varied. Figure 5 depicts boxplots of the time-averages of the four metrics across a range of numbers of observations $I$ and diffusion constants $\alpha$. All combinations of $I \in \{10, 20, 40, 60\}$ and $\alpha \in \{0.001, 0.01, 0.1\}$ were analyzed. These numerical examples confirm our expectations: for advection-dominated dynamics and sparse but low-uncertainty observations, the DLF does significantly better in terms of Mass errors. In terms of the RMS, CoM, and Calibration, we see that the DLF outperforms the KF when data is sparse, namely $I = 10, 20$. The KF gains an advantage when observations are plentiful.

#### 2) Comparison of the KF and DLF Estimates When Uncertainties in the Initial Conditions Are Present

Since the DLF constantly imparts phase information via the likelihood of the pseudo-observations conditioned on the model, the DLF approach should deliver better predictions than the KF on problems where there are uncertainties in the initial data. Again, this is expected when the data has low uncertainty and the dynamics are advection dominated. We will compare the quality metrics of the two methods with initial condition uncertainty. We will show that the DLF manages to overcome this restriction within a few time steps, practically reaching error values comparable to the case of known initial data.

To simulate this uncertainty in the initial data, the truth will be generated with the same initial data as in the previous section, namely $\sigma = 1$ and $\theta = 0.5$, while both the KF and the DLF will be provided different initial data. We will examine two different cases in this section. First, the amplitude $\sigma$ provided to the filters will be drawn uniformly from $\mathcal{U}[\frac{1}{2}, \frac{3}{2}]$, unless individual runs are discussed, for which $\sigma \neq 1$ is picked by hand. Throughout this first set of examples, the phase $\theta = \frac{1}{2}$ is assumed known. In the next set of examples, the phase $\theta$ will be uniformly drawn from $\mathcal{U}[0, 1]$, while assuming $\sigma = 1$ is known. There will again be an exception when discussing individual runs, for which $\theta \neq 0.5$ is picked by hand. The initial data provided to the KF and DLF will be identical for each run to guarantee comparability.

We test the same values for the relative diffusion $\alpha = 0.01$ and the number of data points per observation $I = 20$ as in the previous section. Noise levels remain the same as in the previous section. We will demonstrate that the DLF is significantly more successful in correcting its incomplete knowledge of the initial data than the KF. All trends observed in the previous section, that is, smaller errors and better calibration for the DLF, persist under these settings.

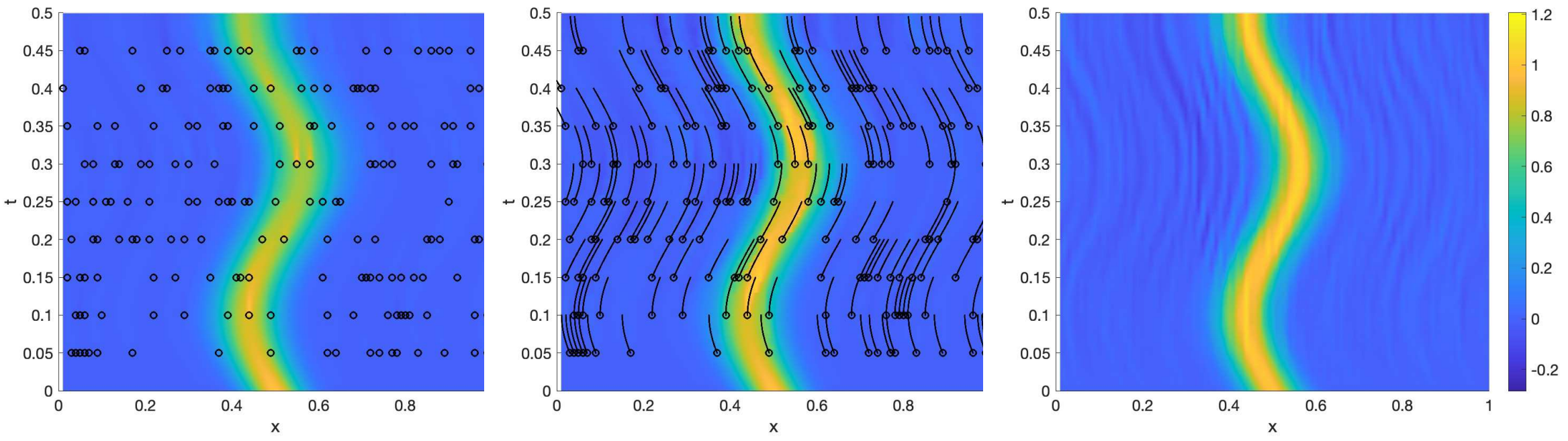

(a) KF (b) DLF (c) Truth

Fig. 3: Posterior mean prediction as estimated by (a) the KF, and (b) the DLF, compared to (c) the truth. Advection dominated case, with $\alpha = 0.01$, initial data $\sigma = 1$ and $\theta = 0.5$ and spatially uncorrelated wave noise $A = 0.05$ and $\tilde{A} = 0$. Both filters use $I = 20$ observations per observation time. The locations of observations are randomly selected grid points, marked by black rings. Observation times are $T_O = \{0.05, 0.1, ..., 0.45\}$. The trajectories of pseudo-observations are shown as black lines.

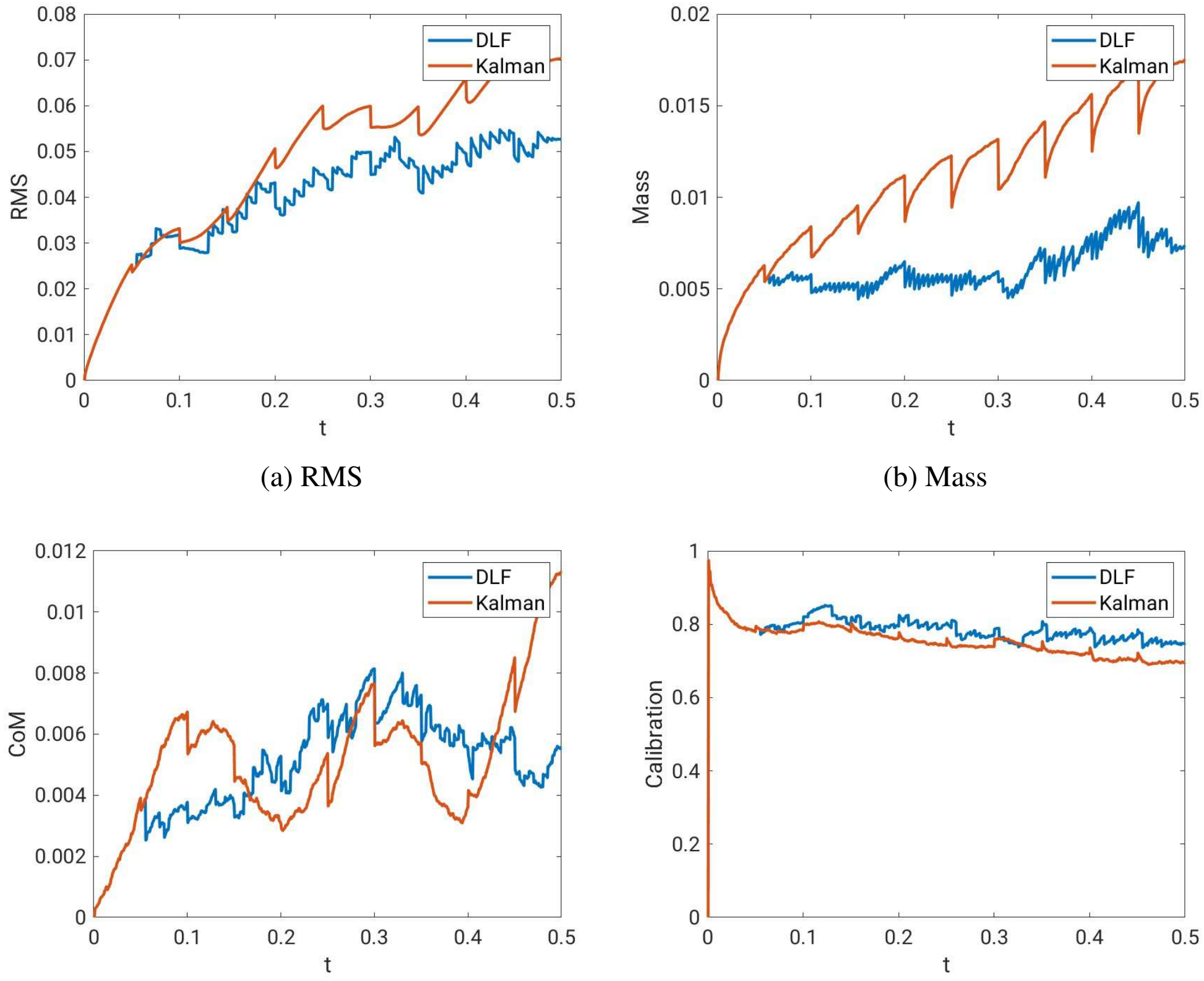

(a) RMS (b) Mass

(c) CoM (d) Calibration

Fig. 4: Time series of (a) RMS, (b) Mass, (c) CoM errors and (d) Calibration of the KF (red), the DLF (blue), averaged over 50 runs. Depicted are results from advection dominated cases with $\alpha = 0.01$, known initial data $\sigma = 1$ and $\theta = 0.5$ and spatially uncorrelated wave noise $A = 0.05$ and $\tilde{A} = 0$. Both filters use $I = 20$ observations per observation time. The locations of observations are randomly selected grid points. Observation times are $T_O = \{0.05, 0.1, ..., 0.45\}$.

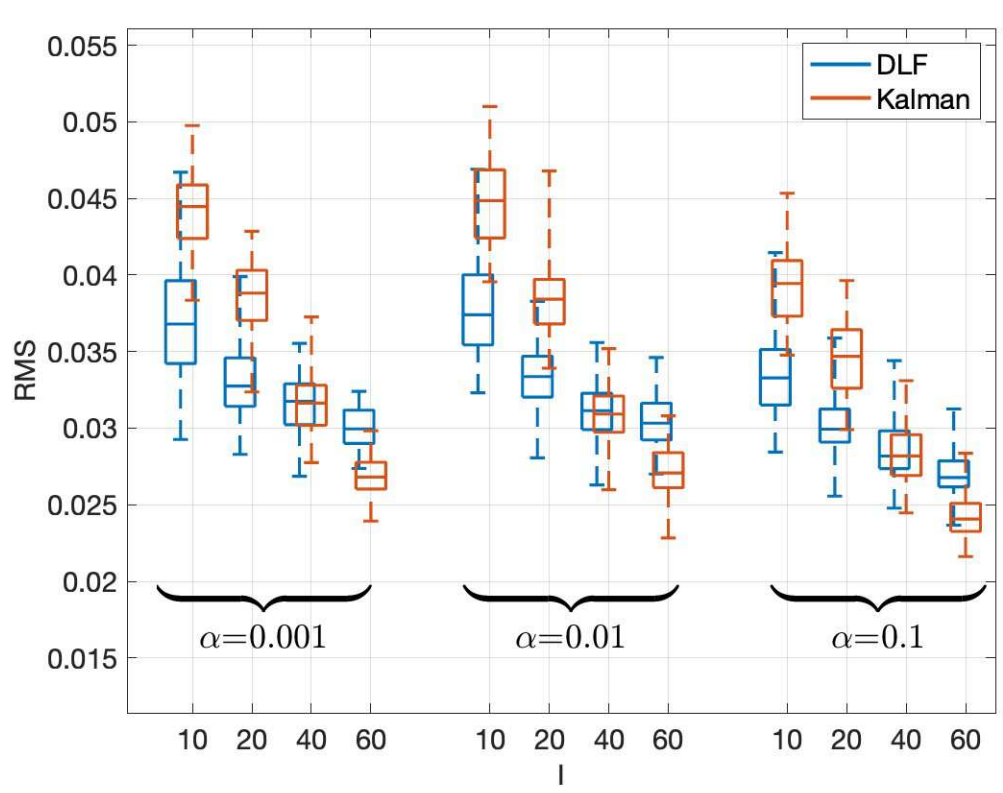

(a) RMS

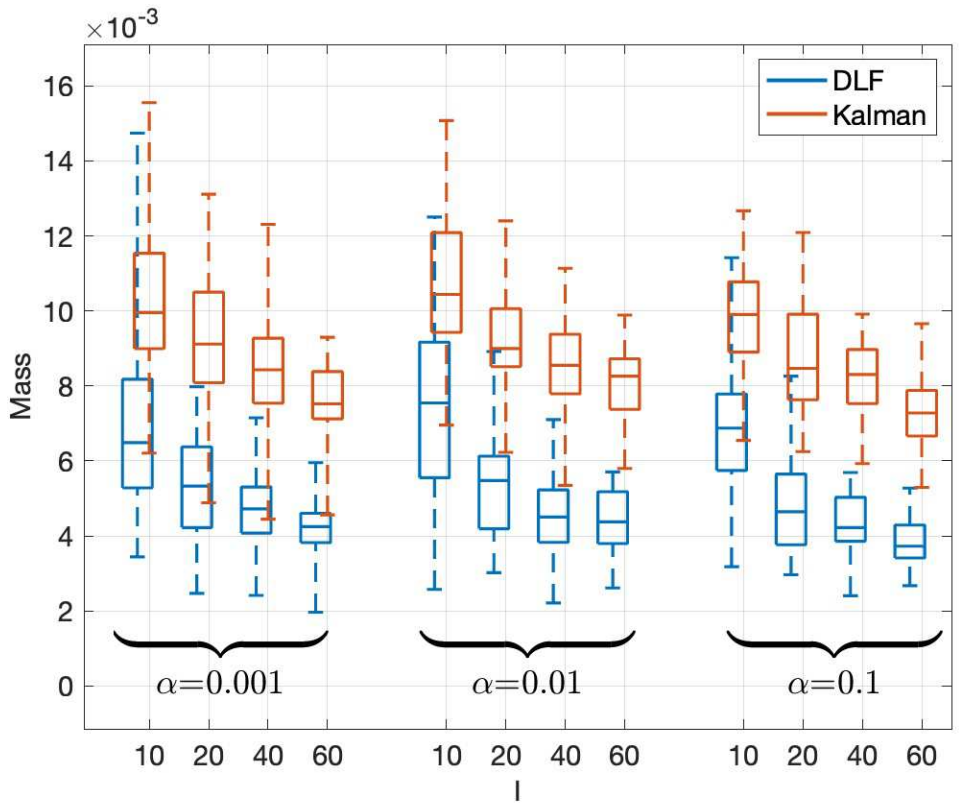

(b) Mass

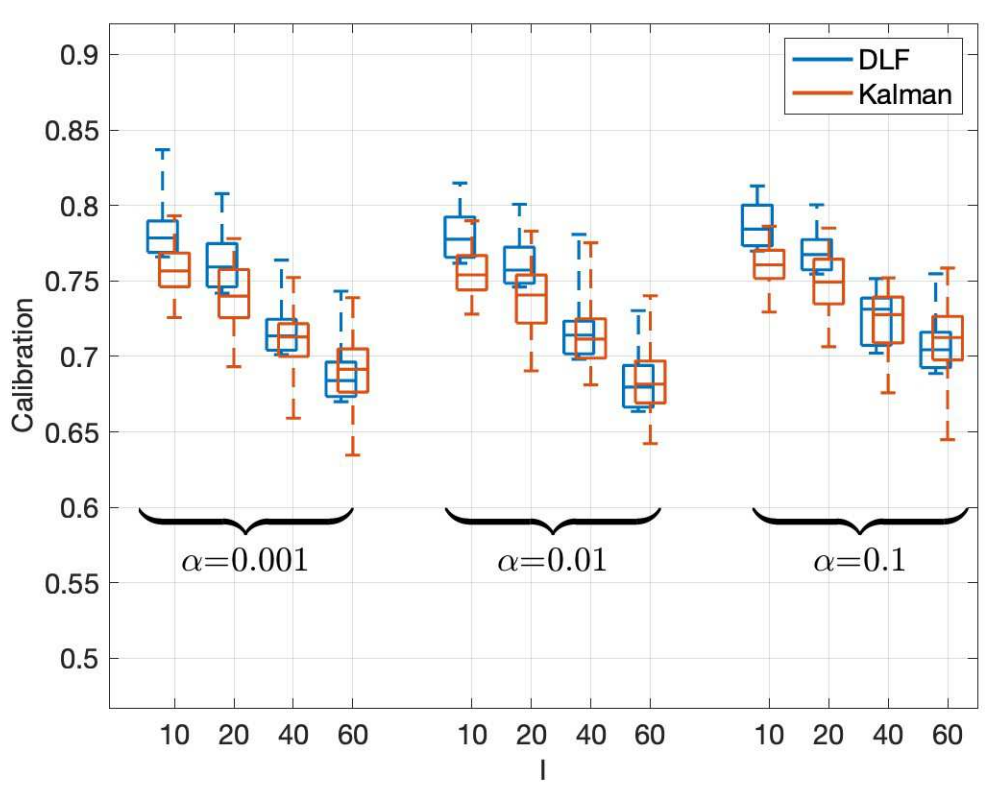

(c) CoM

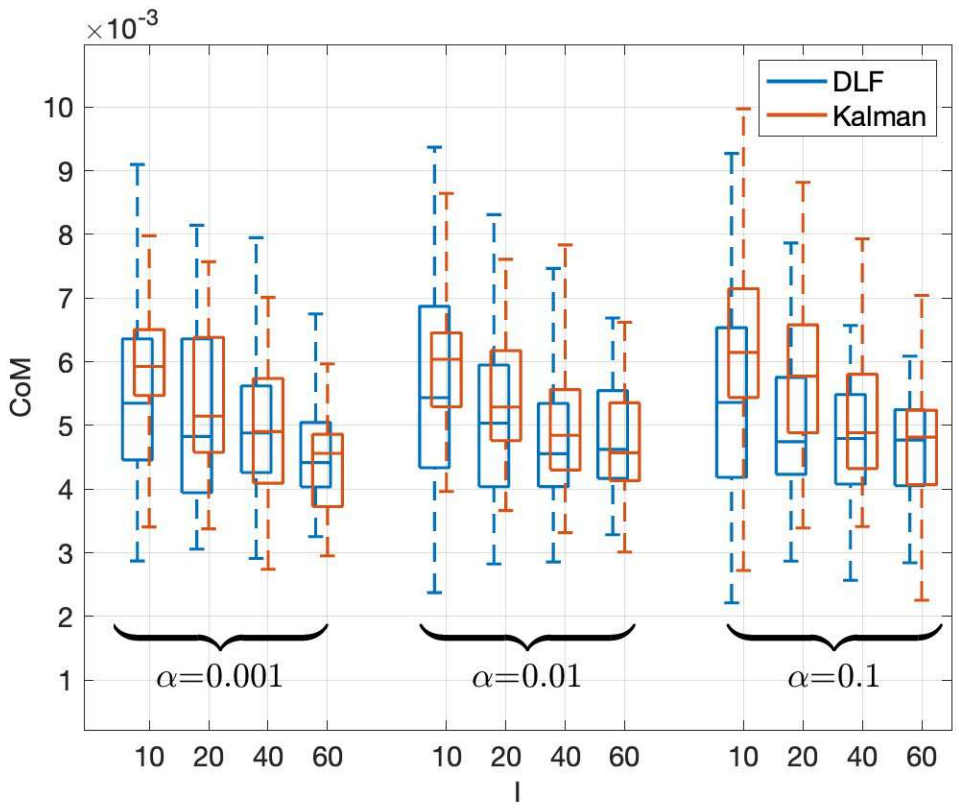

(d) Calibration

Fig. 5: Average (a) RMS, (b) Mass, (c) CoM errors and (d) Calibration of KF (blue), and DLF (red), across 50 runs for spatially independent phase speed noise $A = 0.05, \tilde{A} = 0$, varying diffusion $\alpha =\in \{0, 0.001, 0.01\}$) and observations at $I = 10, 20$ and 40 random locations at every observation time $T_O = \{0.05, 0.1, ..., 0.45\}$. Initial data is assumed known with $\sigma = 1$ and $\theta = \frac{1}{2}$.

In Figure 6, we see the results of an individual run where $\sigma$ was set to 0.7 for the model and $\theta = \frac{1}{2}$. The number of observations available is $I = 20$, and $\alpha = 0.01$.

We see that the DLF manages to correct its incorrect initial amplitude almost immediately as soon as it has access to observations, while the KF struggles to correct its estimation of the amplitude throughout the run.

The time series of the four metrics are shown in Figure 7. These still used $I = 20$ nor $\alpha = 0.01$ was changed from the previous run. We observe qualitatively similar results as before. The DLF performs better in terms of RMS and Mass error, as well as in Calibration. There is no significant difference in the CoM error.

The advantage of the DLF over the KF in terms of RMS and Mass error is significant. After a brief adjustment period, the DLF reaches the same levels as in the case with known initial data. The advantage in Calibration is less pronounced and comparable to the previous case.

Figure 8 shows boxplots of the total value of each of the four metrics across parameters $\alpha \in \{0.001, 0.01, 0.1\}$ and $I \in \{10, 20, 40, 60\}$. We observe the DLF performing better on RMS and Mass error, as well as Calibration, with the advantage expectedly dwindling as the number of observations increases and the diffusion ramps up. In the case of the Mass error, the DLF is still superior across all $I$ and $\alpha$.

Next, we will focus on phase error effects on the estimation. In the remainder of this section, we repeat the previous experiments, but now fix $\sigma = 1$, while $\theta \sim \mathcal{U}[0, 1]$, unless focusing on just a single run.

Figure 9 depicts results from such an individual run with $\theta = 0.25$. The diffusion was, again, set to $\alpha = 0.01$

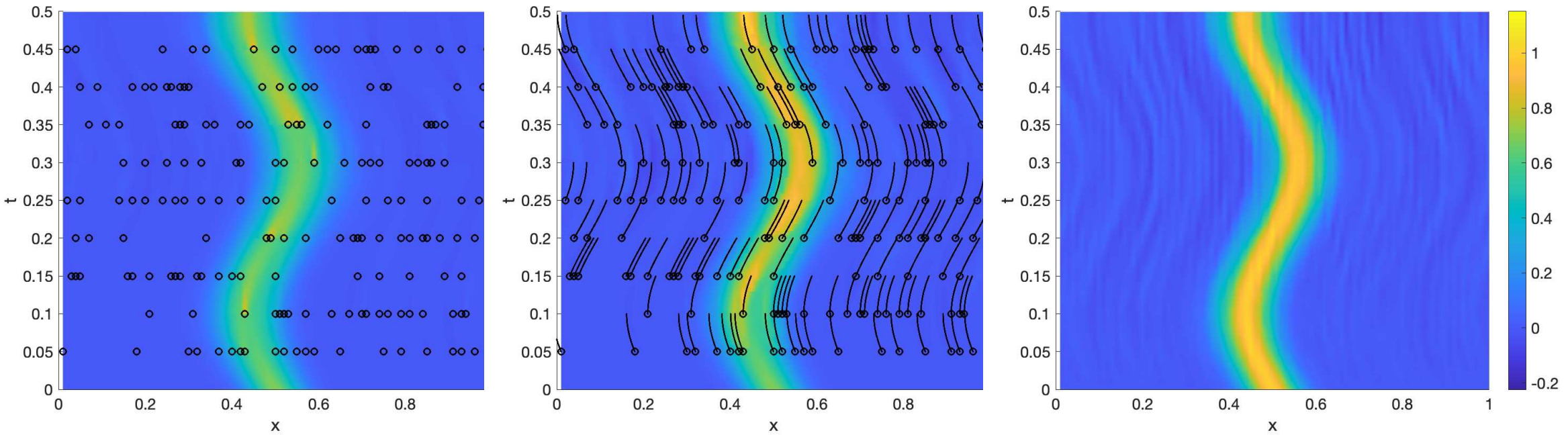

(a) KF (b) DLF (c) Truth

Fig. 6: Posterior mean prediction as estimated by the (a) KF, and the (b) DLF, compared to the (c) truth. Advection dominated case with $\alpha = 0.01$, known initial phase $\theta = 0.5$ and spatially uncorrelated wave noise $A = 0.05$ and $\tilde{A} = 0$. The models use an incorrect initial amplitude of $\sigma = 0.7$ as opposed to the initial amplitude of the truth $\sigma = 1$. Both filters use $I = 20$ observations per observation time. The locations of observations are randomly selected grid points, marked by black rings. Observation times are $T_O = \{0.05, 0.1, \ldots, 0.45\}$. The trajectories of pseudo-observations are shown as black lines.

and $I = 20$. Both filters start with the mode in the wrong position and pick up on their location error as observations become available. The initially displaced mode decreases in amplitude for both models, while a second mode in the correct location starts emerging as soon as observations are available.

The DLF manages to suppress the wrong mode over just a few time steps, roughly the same amount of time it takes to pick up the correct position and amplitude of the actual mode. The KF does significantly worse in this setting, taking almost the entire simulation time to drop the phase error.

In Figure 10, we analyze the average time series of four quality metrics again. The diffusion $\alpha = 0.01$ and number of observations $I = 20$ remain the same. The DLF still performs better than the KF in terms of RMS error and Calibration, with a much bigger advantage in Calibration than in the previously discussed examples. For the first time, there is a clear difference in the CoM error, with the DLF taking a significant lead. Note that in terms of the Mass error, the DLF initially does worse than the KF. This can be explained by the DLF picking up the phase and amplitude of the correct mode, before phasing out the incorrect first mode. In fact there is a brief period where it estimates the existence of two modes. In the long run, the DLF outperforms in terms of Mass error as well.

To close out this section, we test whether the advantage of the DLF over the KF can be sustained for different values of $I$ and $\alpha$. Figure 11 depicts the statistics of the four metrics for $I \in \{10, 20, 40, 60\}$ and $\alpha \in \{0.001, 0.01, 0.1\}$. We observe the following: The DLF does substantially better on RMS, CoM, and Calibration, confirming all trends seen so far in this section. Unlike previously observed, the DLF sustains its advantage throughout cases with higher numbers of observations, likely because more observations make it more probable to pick up the correct location of the true mode in the analysis stage of the assimilation.

In terms of Mass, the DLF now actually performs worse than in previous cases, and even slightly worse than the KF, when few observations are available. This can be explained by the fact that we are looking at time averages of the Mass error here. Since the DLF yields two mode estimates early on, for a brief period, its Mass error is higher. Additional observations help depress the second mode.

### 3) Comparing the KF and DLF When Uncertainties in Phase Speed are Present

In the previous section, we compared the DLF to the KF assuming limited knowledge of the initial data. In this section, we will showcase how both filters perform under the assumption that phase errors are significant. To this end, during the simulation of the truth, $\tilde{A} = 1$ will be used, while the model and thus KF and DLF estimators are unaware of this and still assume $\tilde{A} = 0$. This will cause substantial divergence between the phase speed of the truth and the phase speed used by the model, resulting in significant displacement of the position of the center of mass, if no assimilation happens. Initial data is assumed to be known. The DLF outperforms the KF in terms of RMS and Mass metrics, but it will also be shown that the DLF can correct the displaced center of mass better than the KF. We first take a look at an individual run again. Predictions and truth are depicted in Figure 12. Diffusion is again $\alpha = 0.01$ and $I = 20$ observations are available at each observation time for this example. As seen in previous sections the DLF manages to adjust its predictions to the observations much more rapidly than the KF.

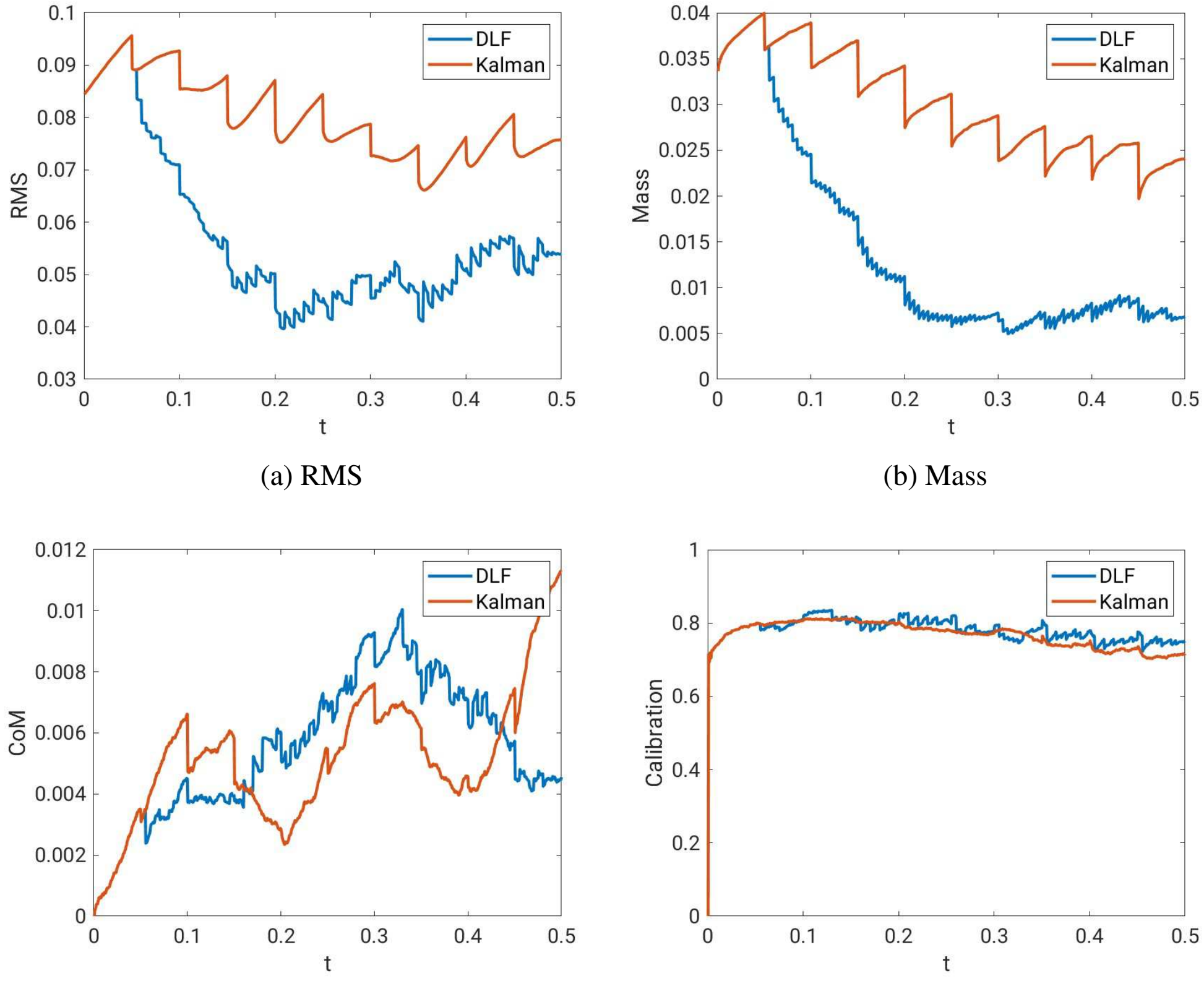

(a) RMS

(b) Mass

(c) CoM

(d) Calibration

Fig. 7: Time series of (a) RMS, (b) Mass, (c) CoM errors and (d) Calibration of the KF (red), the DLF (blue), averaged over 50 runs. Depicted are results from advection dominated cases with $\alpha = 0.01$, known initial phase $\theta = 0.5$ and spatially uncorrelated wave noise $A = 0.05$ and $\tilde{A} = 0$. The models use an incorrect initial amplitude of $\sigma = 0.7$ as opposed to the initial amplitude of the truth $\sigma = 1$. Both filters use $I = 20$ observations per observation time. The locations of observations are randomly selected grid points. Observation times are $T_O = \{0.05, 0.1, ..., 0.45\}$.

Regarding the time series of the quality metrics depicted in Figure 13 we now see the DLF outperforming the traditional KF in all metrics except the CoM error, where it occasionally does slightly worse.

Considering the total values of the four metrics over a larger range of $\alpha$ and $I$ we see the DLF clearly taking the lead. Boxplots of these statistics are depicted in figure 14. The DLF performs better on average regarding all four metrics, this time persisting through higher diffusion and increased number of observations, as far as examined. This again shows the superiority of the DLF over the KF in the case of an ill-informed model, this time illustrated by the models getting the phase speed wrong.

#### 4) Dynamics and the DLF and KF outcomes

In previous sections, we compared the capabilities of the DLF and the KF under increasingly complex uncertainties in the model. All these experiments were initially conducted in the advection-dominated setting $\alpha \ll 1$. In this next section, we will consider the performance of the DLF and KF under different dynamic conditions, i.e., when $\alpha \approx 1$ or larger. The DLF relies on the propagation of observations along characteristics determined by the advective part of the system. The evolution of these observations along these characteristics occurs in the presence of diffusion, which is determined by the derivatives of the imperfect model $v$. Thus, we expect diminishing returns in the high $\alpha$ regime for the DLF. We will examine this next. In the following experiments, $I = 20$, $\sigma = 1$, $\theta = \frac{1}{2}$, and the noise remains spatially uncorrelated, i.e., $\tilde{A} = 0$. The resulting average metrics as a function of $\alpha$ on a range of $\alpha \in [0.0001, 5]$ are depicted in Figure 15. As noted previously, the DLF has an advantage over the KF in all four analyzed metrics as long as advection dominates, i.e., $\alpha \ll 1$. However, this advantage decreases as $\alpha$ increases. Nonetheless, the DLF remains useful if the data is sparse.

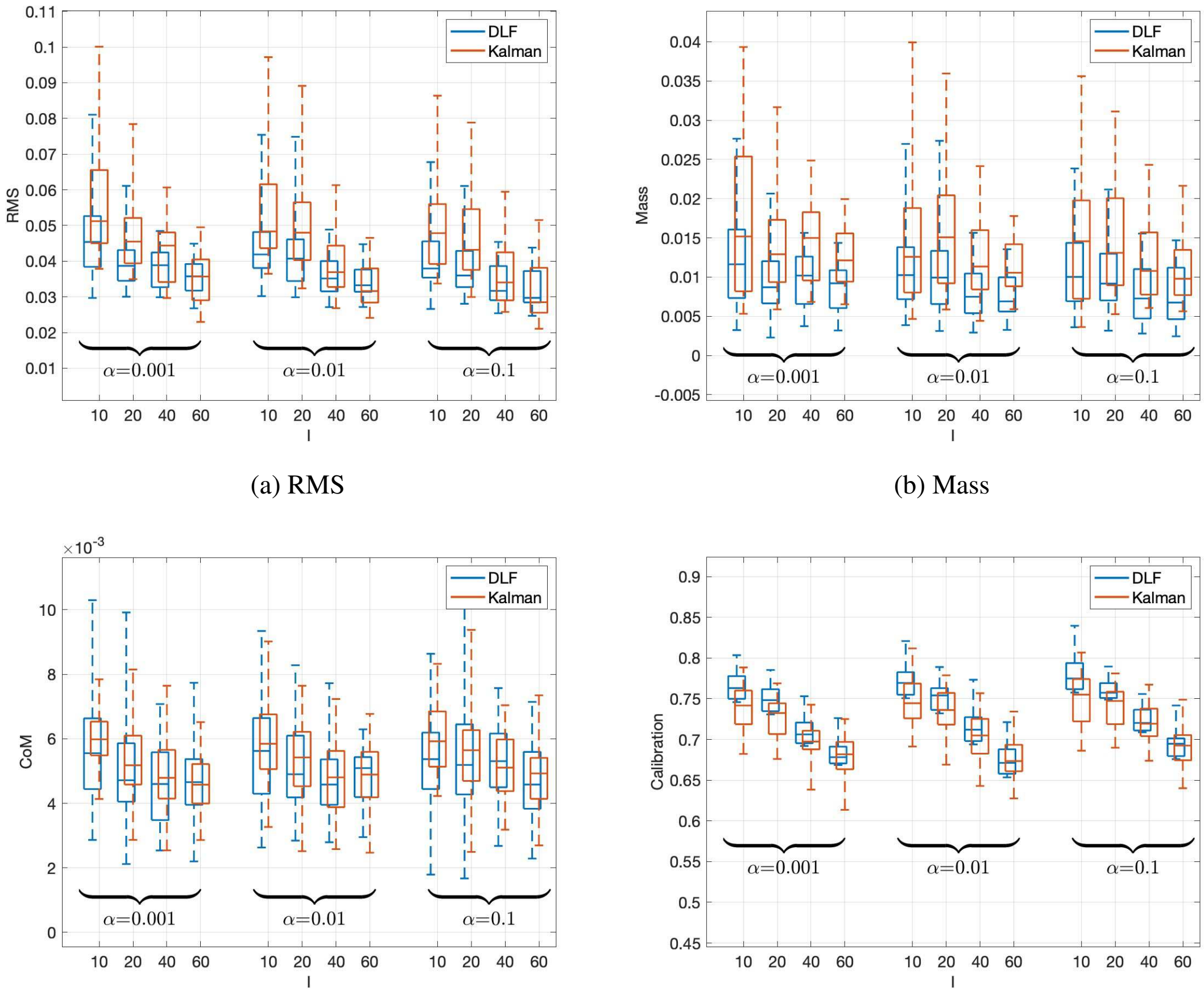

(a) RMS

(b) Mass

(c) CoM

(d) Calibration

Fig. 8: Average (a) RMS, (b) Mass, (c) CoM errors and (d) Calibration of KF (blue), DLF (red), across 50 runs for spatially independent phase speed noise $A = 0.05, \tilde{A} = 0$, varying difusion $\alpha =\in \{0, 0.001, 0.01\}$) and observations at $I = 10, 20$ and 40 random locations at every observation time $T_O = \{0.05, 0.1, ..., 0.45\}$. Initial amplitude is assumed uncertain with $\sigma = \mathcal{U}\left(\frac{1}{2}, \frac{3}{2}\right)$ and $\theta = \frac{1}{2}$.

In the advection-dominated case, the advantages of the DLF were more pronounced in cases of ill-informed models. Thus, as a last experiment, we will investigate if these advantages can be maintained as $\alpha$ increases. To this end, we assume uncertainties in the initial data amplitude and phase, i.e., $\sigma \sim \mathcal{U}[\frac{1}{2}, \frac{3}{2}]$ and $\theta \sim \mathcal{U}[0, 1]$. Further, we reintroduce the systematic uncertainty in phase speed into the model, i.e., $\tilde{A} = 1$, when simulating the truth. The number of observations is set to $I = 20$, while $\alpha \in [0.0001, 5]$. The resulting average metrics are depicted in Figure 16. We see now that the DLF, again, performs better in all four metrics over the entire range of analyzed $\alpha$. For RMS and Mass error, as well as Calibration, the distance between DLF and KF decreases as $\alpha$ increases, while the advantage in terms of CoM error seems to be nearly unaffected by the values of $\alpha$.

## 6. Discussion and Conclusions

The DLF approach to data assimilation was developed to handle hyperbolic (wave) dynamics. In this work, we extend the DLF approach to handle advection-diffusion dynamics. The DLF approach was first proposed in Restrepo (2013) and made operational in Foster and Restrepo (2022) on hyperbolic problems. A significant challenge in extending DLF to advection-diffusion problems is that the diffusion term needs to be evaluated along characteristic paths. We used the derivatives of the model as an estimator for this term. However, there are other alternatives, depending on the physics underlying the problem.

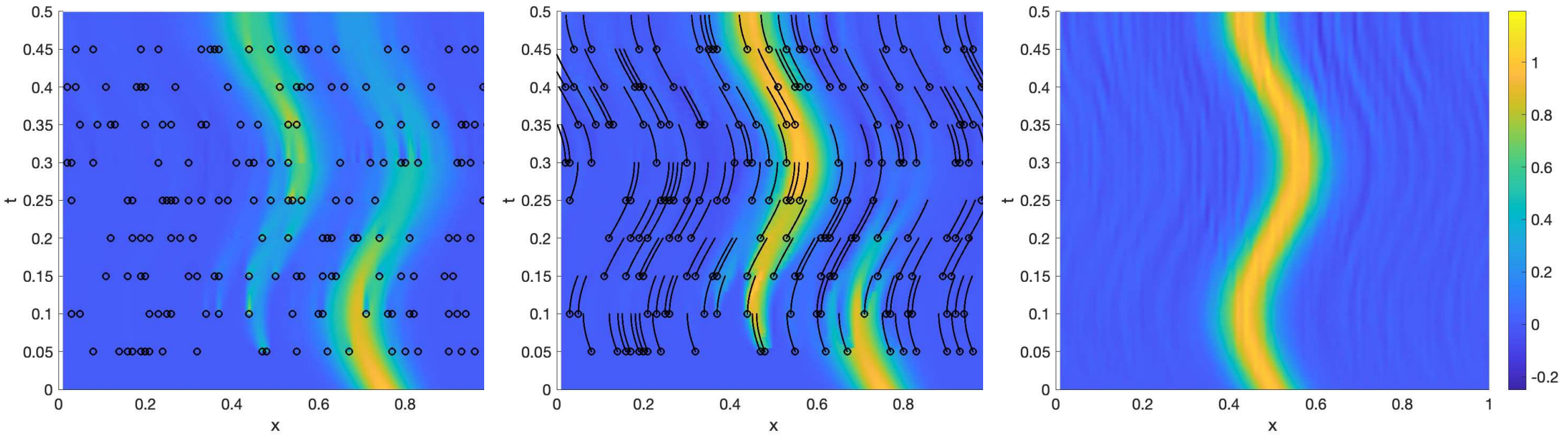

(a) KF (b) DLF (c) Truth

Fig. 9: Posterior mean prediction as estimated by the (a) KF, and the (b) DLF, compared to the (c) truth. Advection dominated case with $\alpha = 0.01$, known initial amplitude $\sigma = 1$ and spatially uncorrelated wave noise $A = 0.05$ and $\tilde{A} = 0$. The models use an incorrect initial phase of $\theta = 0.25$ as opposed to the initial amplitude of the truth $\theta = 0.5$. Both filters use $I = 20$ observations per observation time. The locations of observations are randomly selected grid points, marked by black rings. Observation times are $T_O = \{0.05, 0.1, ..., 0.45\}$. The trajectories of pseudo-observations are shown as black lines.

In this study, the evolution of observations through time was driven by a combination of data and the model at hand. Since the grid points used in model calculations and the characteristics along which observations travel in the Lagrangian frame did not necessarily coincide, extrapolation of the model to these off-grid points was required.

There are two ways the DLF can be formulated. In both cases, phase information is conveyed via the dynamic likelihood. In the multi-analysis case, observations or their projections forward in time, along with their uncertainties, can improve the extent of space in which observations affect the analysis product. In practice, the decision of which observations to keep for how long after their original measurement will need to balance computational complexity against improved accuracy.

The DLF approach requires a code that can solve the characteristics problem. With this solver, all additional implementation steps are no harder to implement than the classic KF. Like the KF, the computational complexity of the DLF approach applied to the KF is cubic in the maximum number of data points used at a time step. It can, however, still be significantly higher than the KF, since the DLF would be potentially applied more frequently than if applied at observation times exclusively. Some countermeasures to keep the DLF's complexity in check were discussed. Improved estimates make the higher computational cost justifiable.

Using numerical simulations, we demonstrated that the dynamic likelihood filter (DLF) outperformed the standard KF estimation concerning several metrics of accuracy. We showed that the DLF is superior to the KF when advection dominates diffusion, and observations are sparse and have high precision. Further, we demonstrated that:

- The DLF leads to a more accurate prediction of the truth than the KF, as demonstrated through its lower RMS in all experiments.
- The DLF estimates are significantly less sensitive to uncertainties in the initial data than the KF. It manages to predict the correct phase and amplitude within a shorter time and does so more accurately. As a result of capturing the phase more accurately, the center of mass of the solution is predicted with more accuracy. Further, the DLF gives more accurate local estimates (RMS) and predictions of overall mass.
- The DLF leads to more accurate predictions even when the phase speed is affected by a great deal of uncertainty.
- The DLF is superior in estimating the variance, particularly when uncertainty beyond noise is introduced through uncertain initial data or an ill-informed model.

The last three points remain true, even when diffusion and advection are roughly the same ($\alpha \approx 1$).

In summary, the DLF approach to data assimilation on advection-diffusion problems shows great promise as an estimator, particularly when the observation network is sparse yet of low noise. The implementation requires special time integrators, but its computational overhead is well offset by producing better estimates. In Foster and Restrepo (2022), we showed that the DLF permits Bayesian estimates of model and pseudo-observations *into the future*, possibly beyond the present time when no observations are available. Using conventional data assimilation, forecasts will use the model-informed prior in the estimate of future moments. If the pseudo-observations inform a likelihood that is more compact than the prior, the forecast of the mean of the state may well be significantly different than the mean predicted via the prior only. This unique capability of being able to make Bayesian forecasts persists in the DLF approach, as applied to advection-diffusion problems.

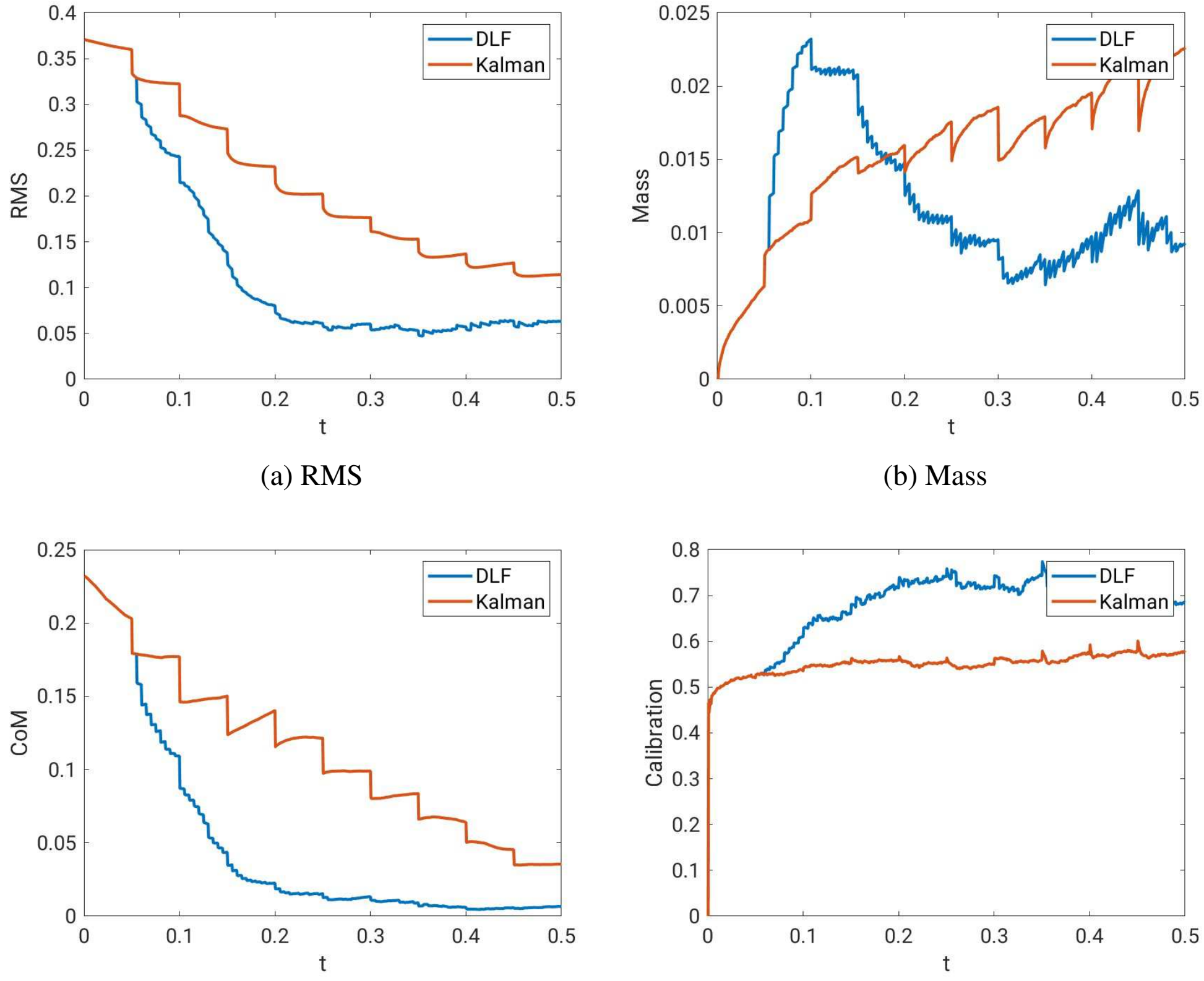

(a) RMS

(b) Mass

(c) CoM

(d) Calibration

Fig. 10: Time series of (a) RMS, (b) Mass, (c) CoM errors and (d) Calibration of the KF (red), the DLF (blue), averaged over 50 runs. Advection dominated cases with $\alpha = 0.01$, known initial amplitude $\sigma = 1$ and spatially uncorrelated wave noise $A = 0.05$ and $\tilde{A} = 0$. The models use an incorrect initial phase of $\theta = 0.25$ as opposed to the initial amplitude of the truth $\theta = 0.5$. Both filters use $I = 20$ observations per observation time. The locations of observations are randomly selected grid points. Observation times are $T_O = \{0.05, 0.1, ..., 0.45\}$.

## 7. Acknowledgments

The submitted manuscript has been authored by a contractor of the U.S. Government under Contract No. DE-AC05-00OR22725. Accordingly, the U.S. Government retains a non-exclusive, royalty-free license to publish or reproduce the published form of this contribution, or allow others to do so, for U.S. Government purposes.

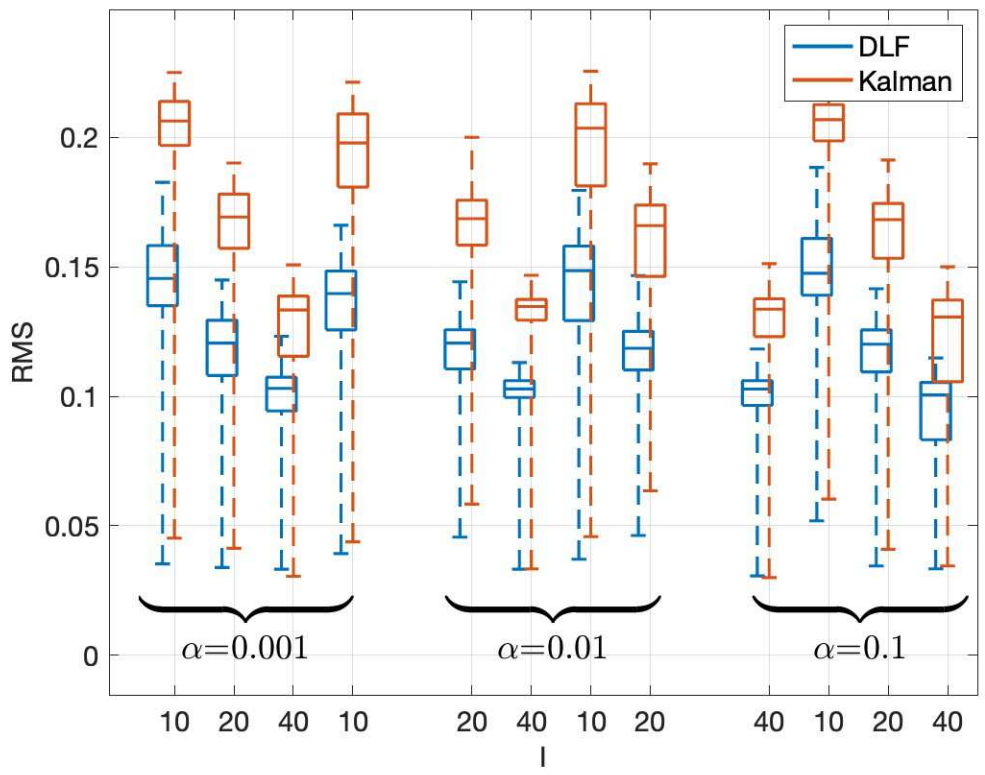

(a) RMS

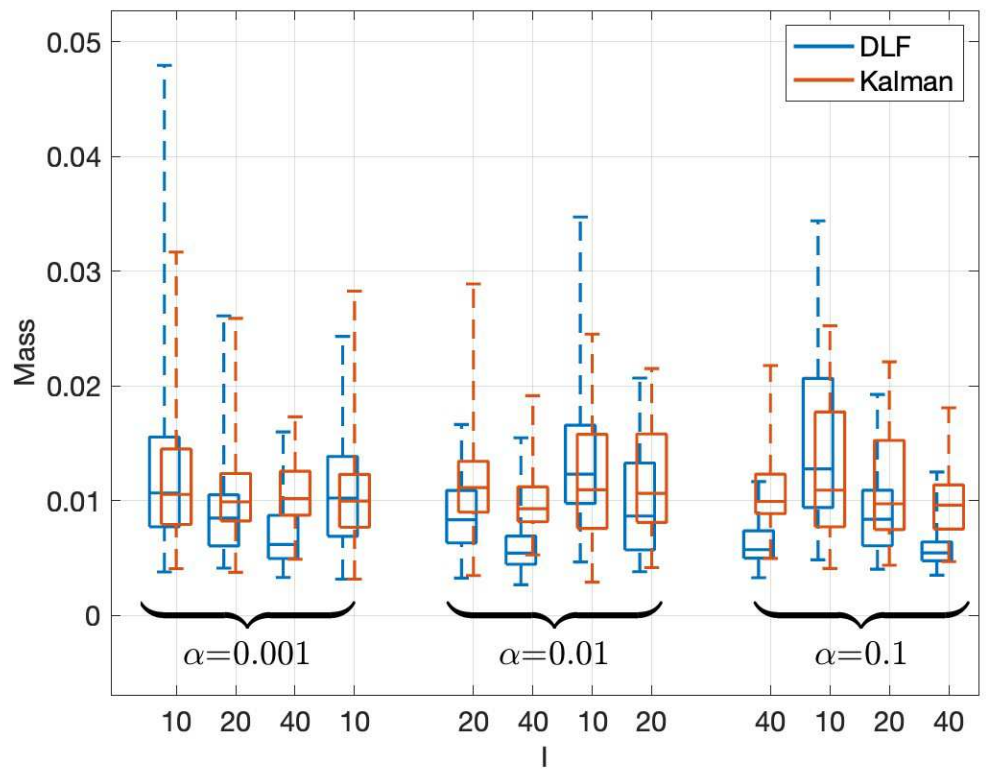

(b) Mass

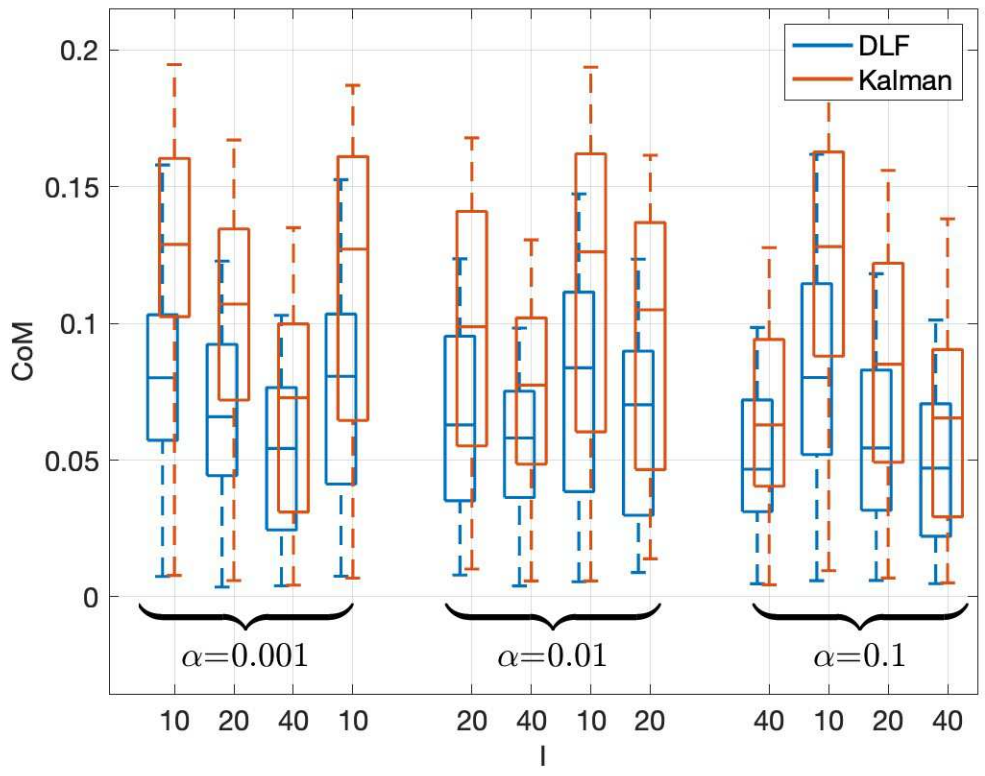

(c) CoM

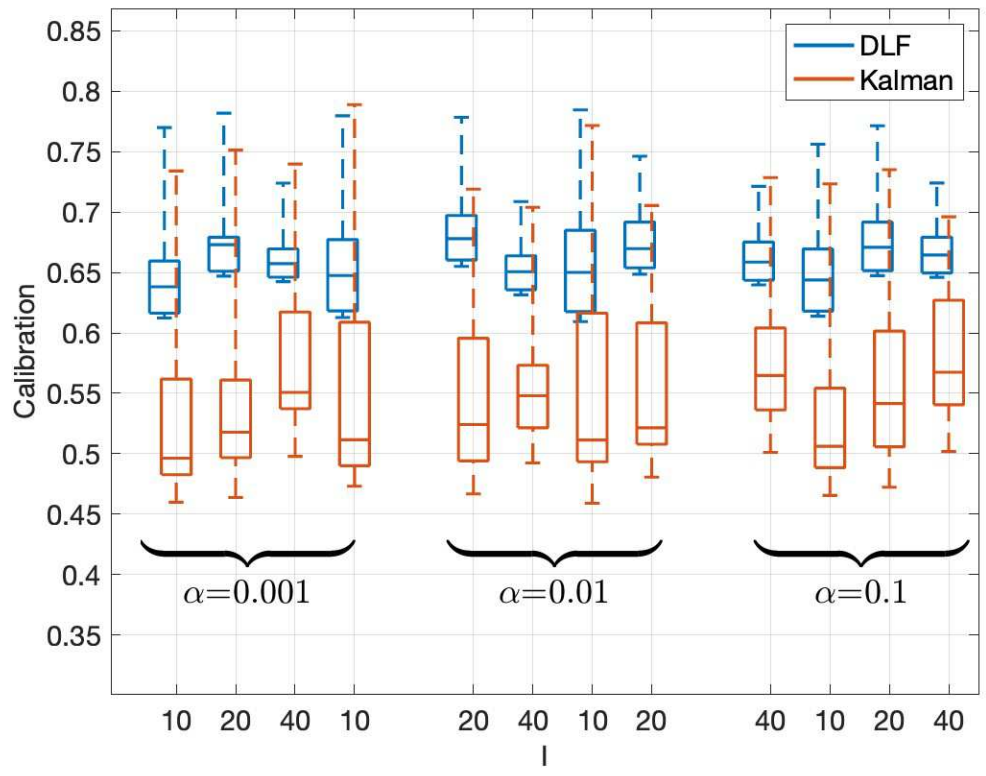

(d) Calibration

FIG. 11: Average (a) RMS, (b) Mass, (c) CoM errors and (d) Calibration of KF (blue), DLF (red), across 50 runs for spatially independent phase speed noise $A = 0.05, \tilde{A} = 0$, varying difusion $\alpha = \in \{0, 0.001, 0.01\}$) and observations at $I = 10, 20$ and 40 random locations at every observation time $T_O = \{0.05, 0.1, ..., 0.45\}$. Initial phase is assumed uncertain with $\sigma = 1$ and $\theta = \mathcal{U}(0, 1)$.

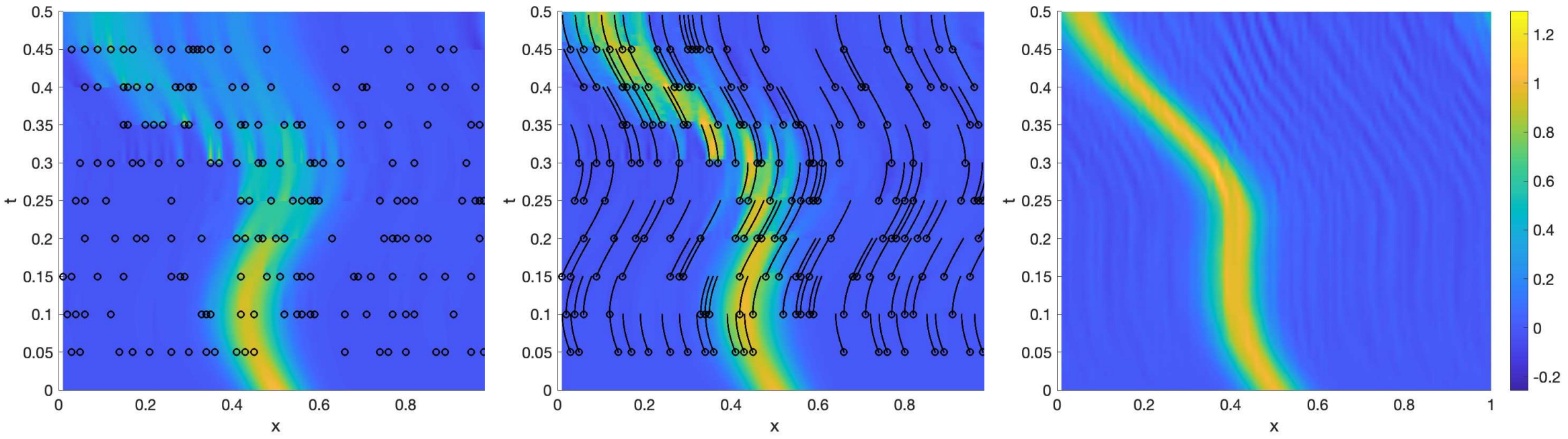

(a) KF (b) DLF (c) Truth

Fig. 12: Posterior mean prediction as estimated by the (a) KF, and the (b) DLF, compared to (c) the truth. Advection dominated case with $\alpha = 0.01$, known initial data $\sigma = 1$ and $\theta = 0.5$. Phase speed noise of the truth is assumed spatially correlated $\tilde{A} = 1$, while both models assume $\tilde{A} = 0$, causing significant discrepancies in phase speed between truth and model. Both filters use $I = 20$ observations per observation time. The locations of observations are randomly selected grid points, marked by black rings. Observation times are $T_O = \{0.05, 0.1, ..., 0.45\}$. The trajectories of pseudo-observations are shown as black lines.

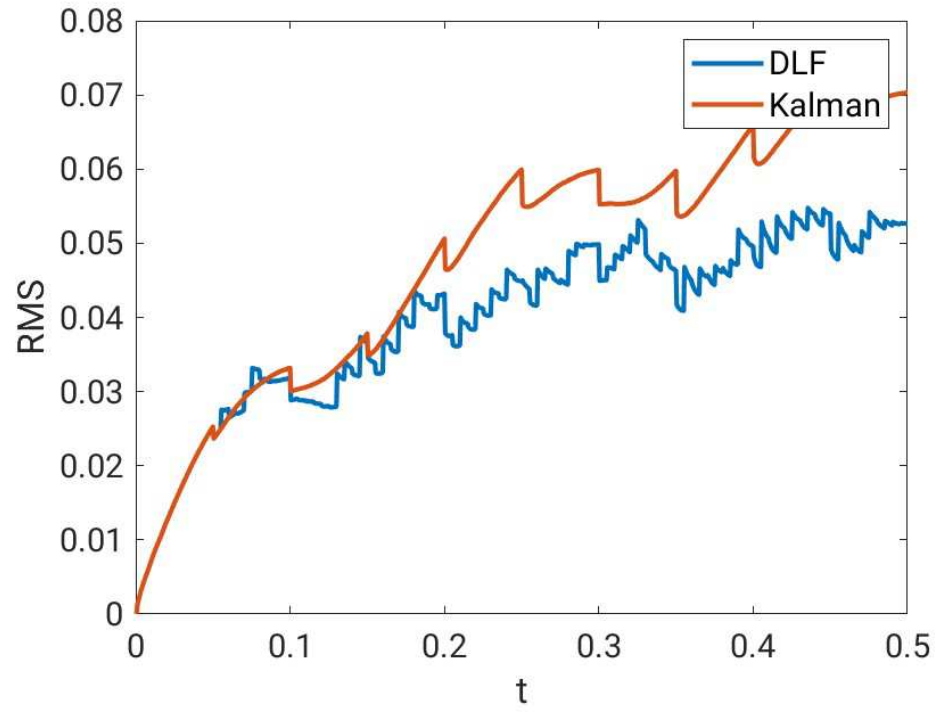

(a) RMS

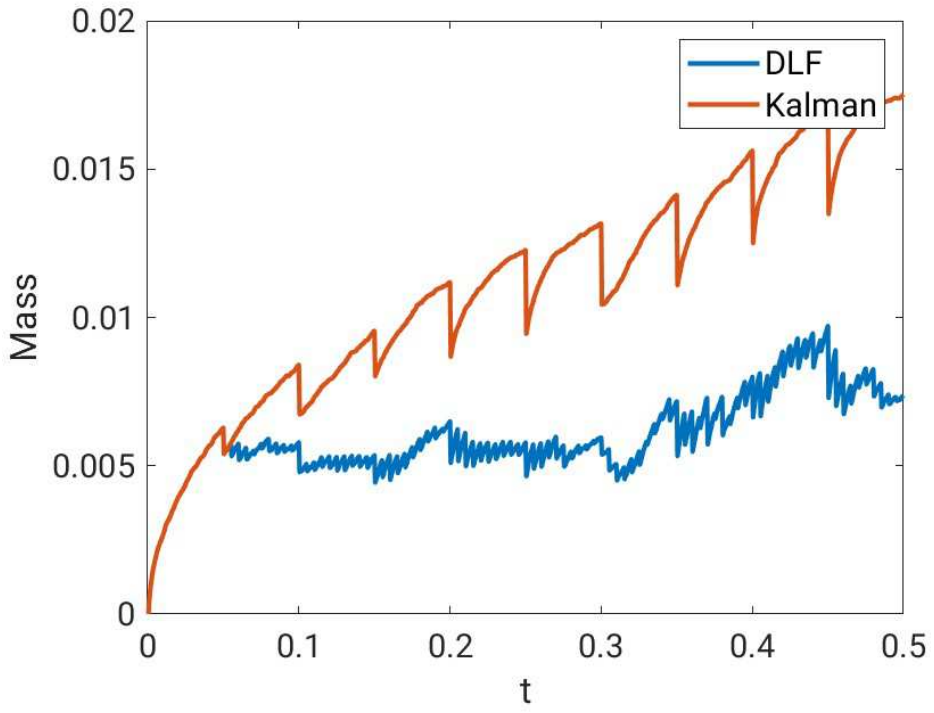

(b) Mass

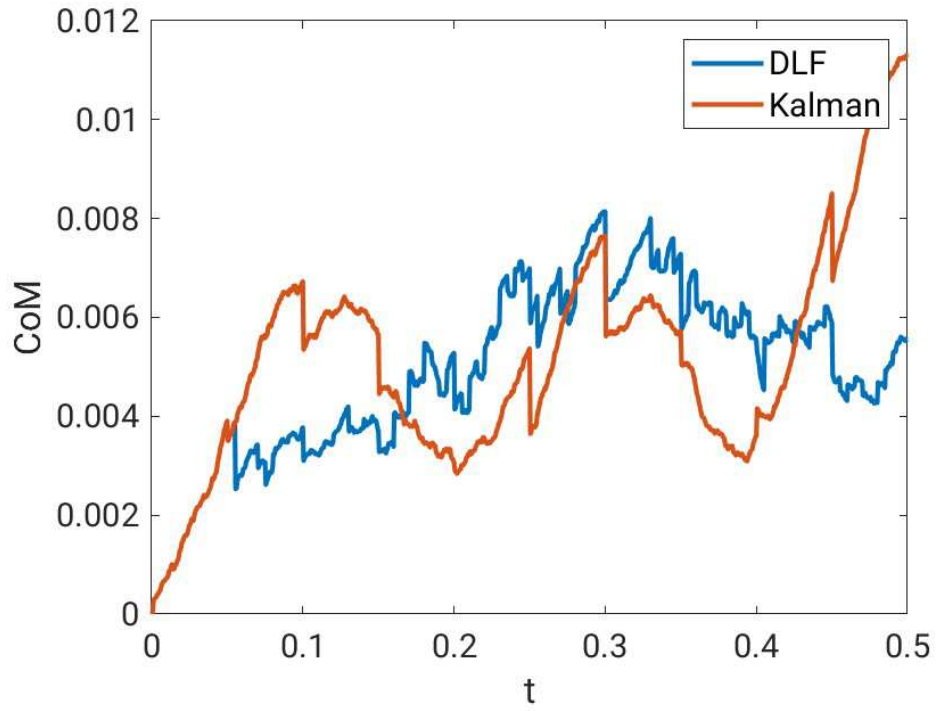

(c) CoM

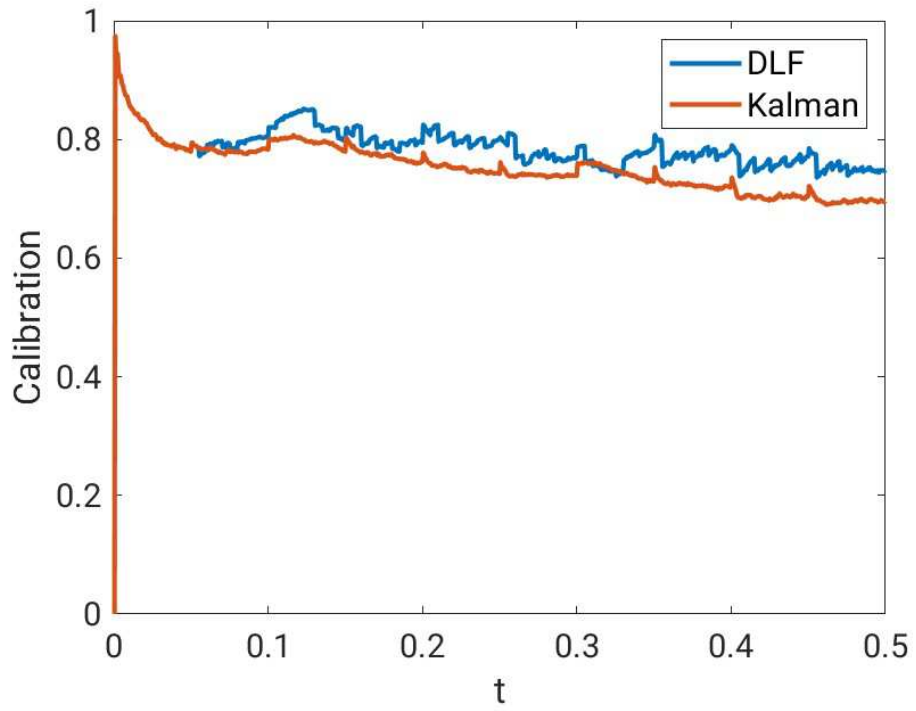

(d) Calibration

Fig. 13: Posterior mean prediction as estimated by the (a) KF, the (b) DLF, compared to (c) the truth. Advection dominated case with $\alpha = 0.01$, known initial data $\sigma = 1$ and $\theta = 0.5$. Phase speed noise of the truth is assumed spatially correlated $\tilde{A} = 1$, while both models assume $\tilde{A} = 0$, causing significant discrepancies in phase speed between truth and model. Both filters use $I = 20$ observations per observation time. The locations of observations are randomly selected grid points. Observation times are $T_O = \{0.05, 0.1, ..., 0.45\}$.

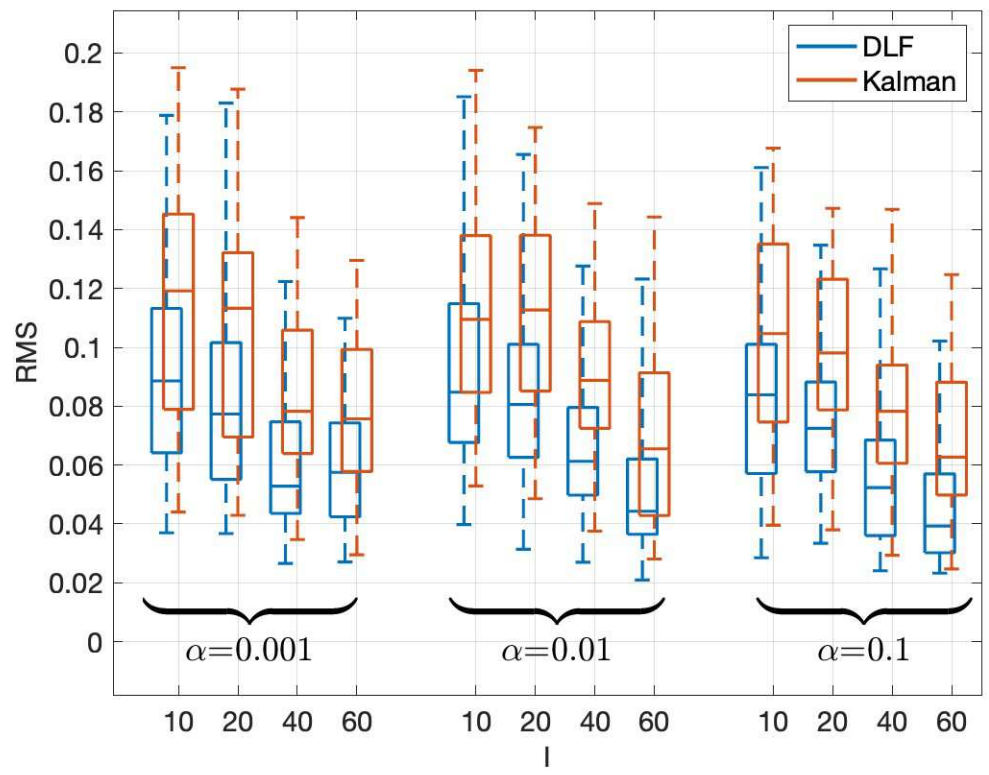

(a) RMS

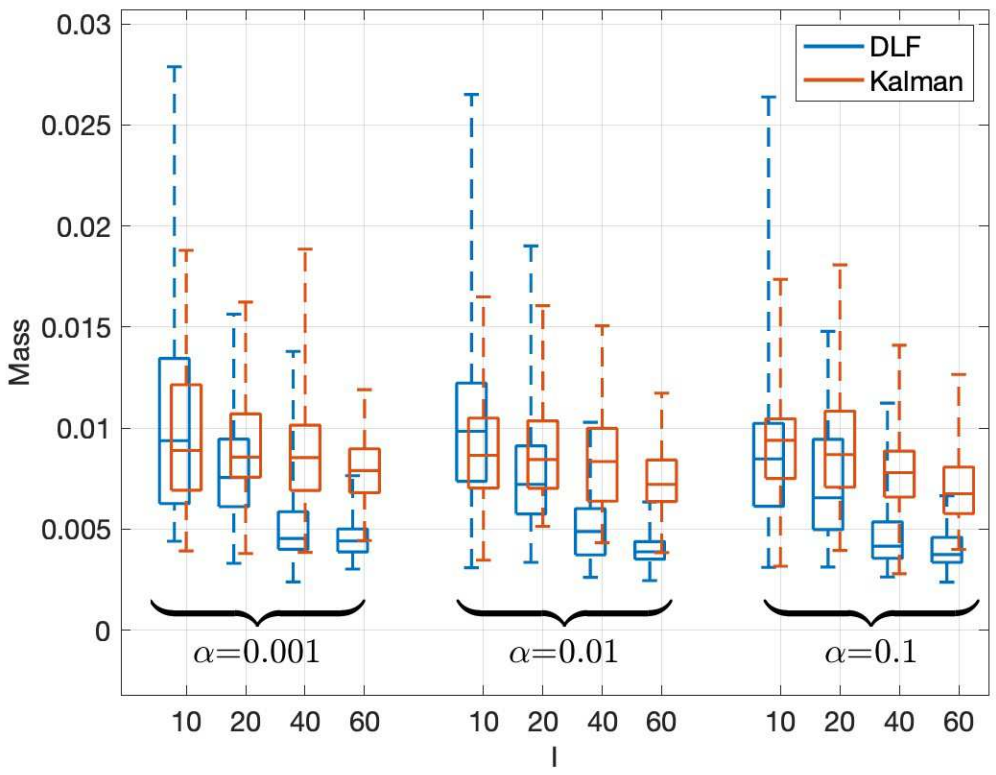

(b) Mass

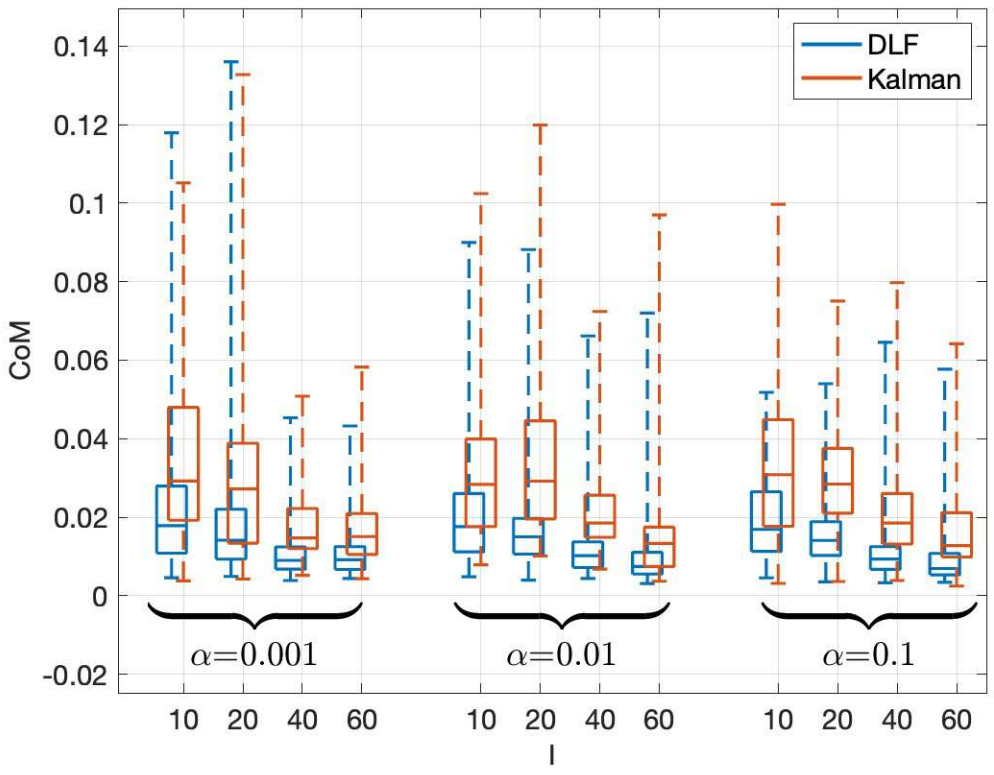

(c) CoM

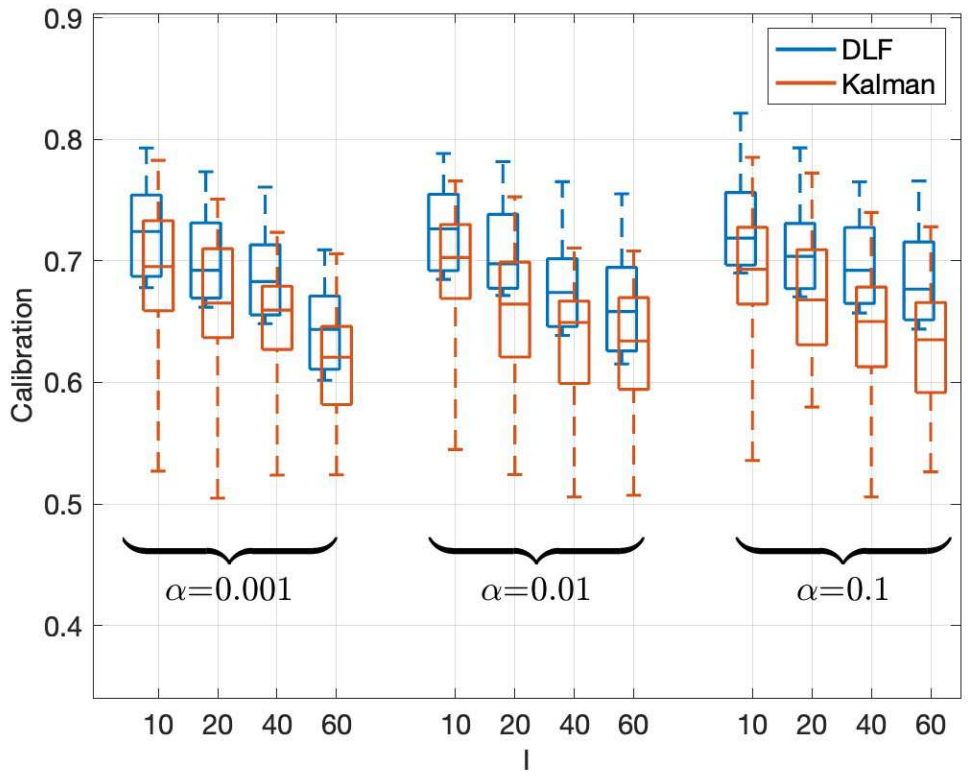

(d) Calibration

Fig. 14: Average (a) RMS, (b) Mass, (c) CoM errors and (d) Calibration of KF (blue), DLF (red), across 50 runs for spatially correlated phase speed noise $A = 0.05, \tilde{A} = 1$, varying difusion $\alpha = \in \{0, 0.001, 0.01\}$) and observations at $I = 10, 20$ and 40 random locations at every observation time $T_O = \{0.05, 0.1, ..., 0.45\}$. Initial data is assumed known initial $\sigma = 1$ and $\theta = \frac{1}{2}$.

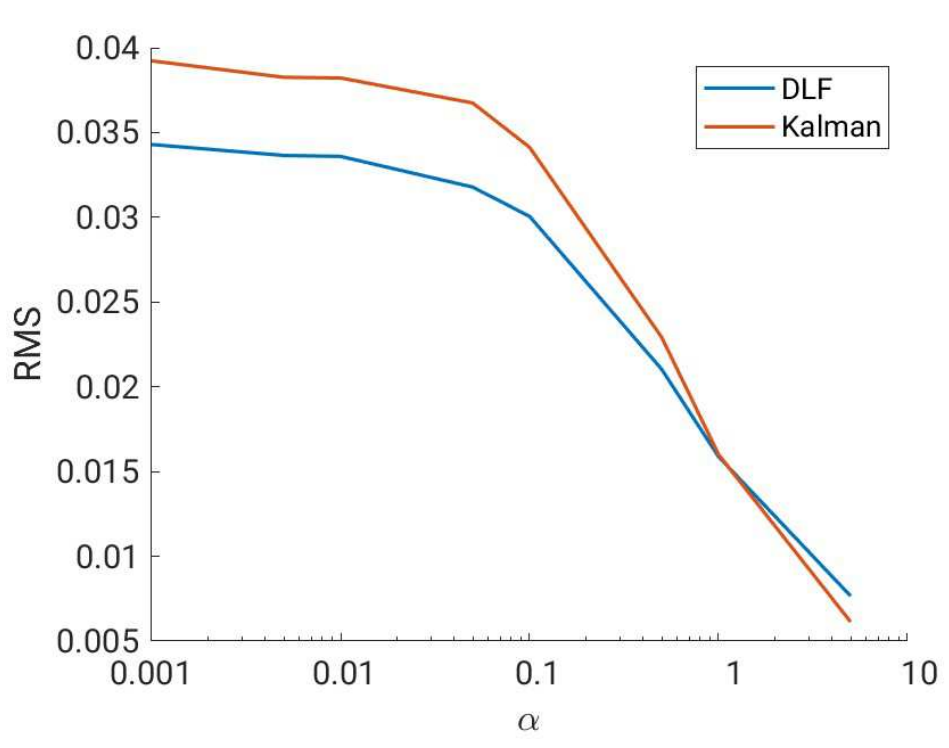

(a) RMS

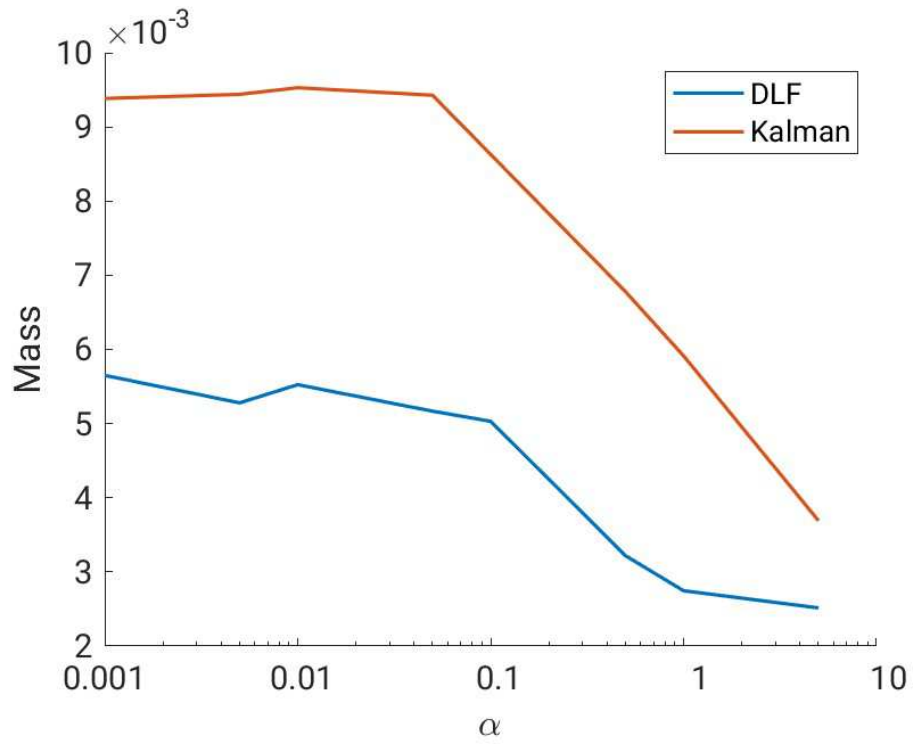

(b) Mass

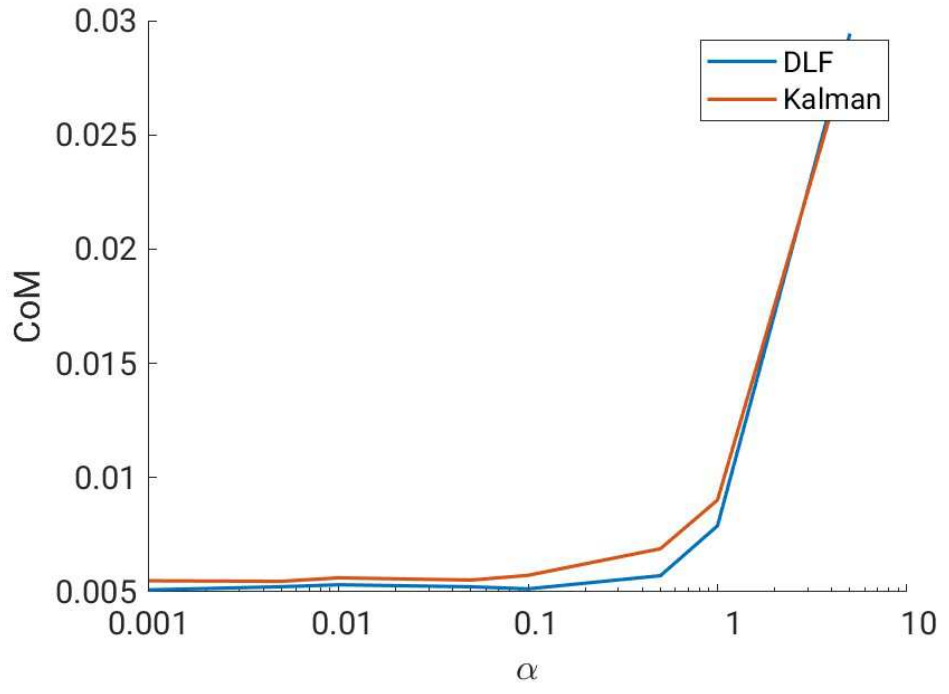

(c) CoM

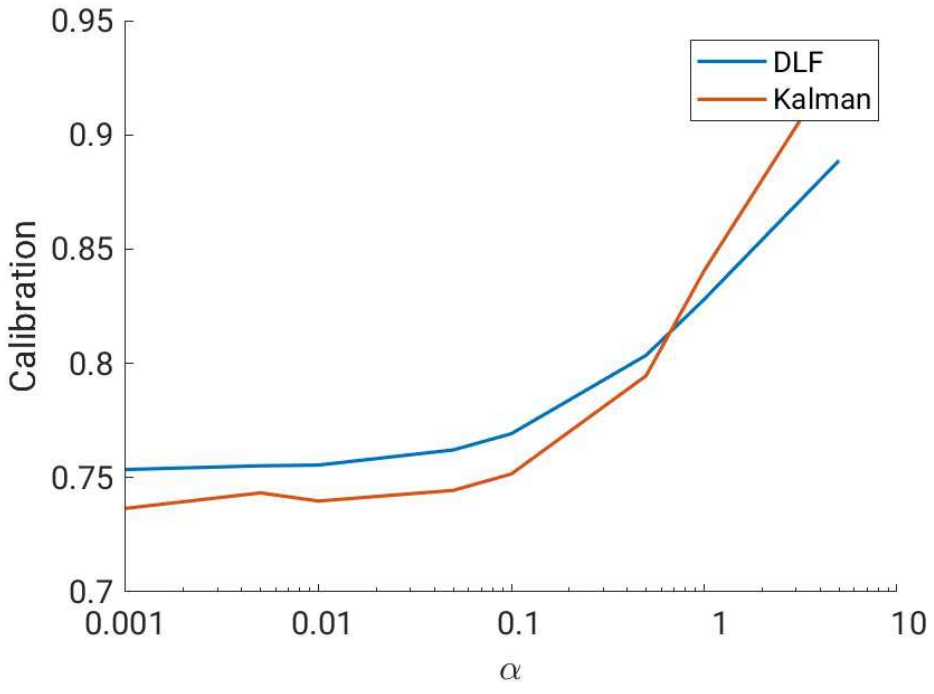

(d) Calibration

Fig. 15: Average (a) RMS, (b) Mass, (c) CoM errors, and (d) Calibration of KF (blue) and DLF (red) as a function of $\alpha$, across 50 runs for spatially independent phase speed noise $A = 0.05$, $\tilde{A} = 0$, and $I = 20$ randomly located observations available at every observation time $T_O = \{0.05, 0.1, ..., 0.45\}$ with known initial data $\sigma = 1$ and $\theta = \frac{1}{2}$.

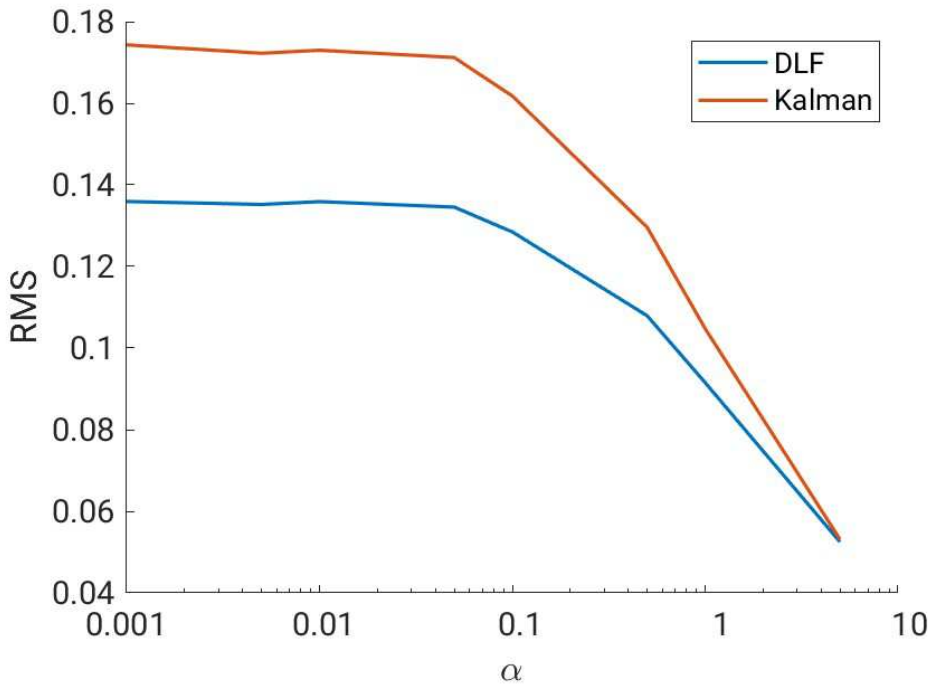

(a) RMS

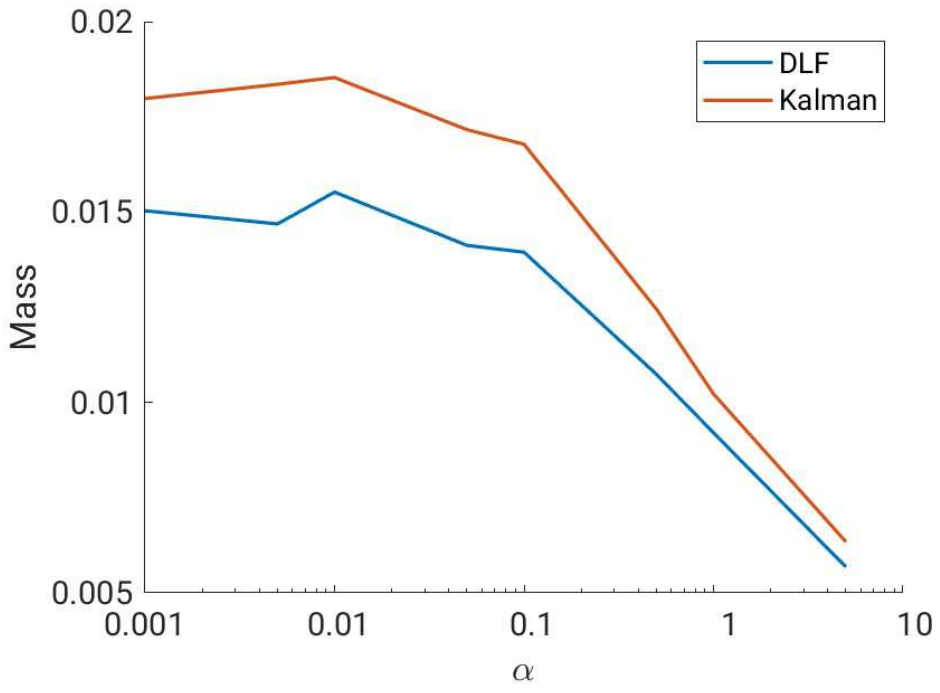

(b) Mass

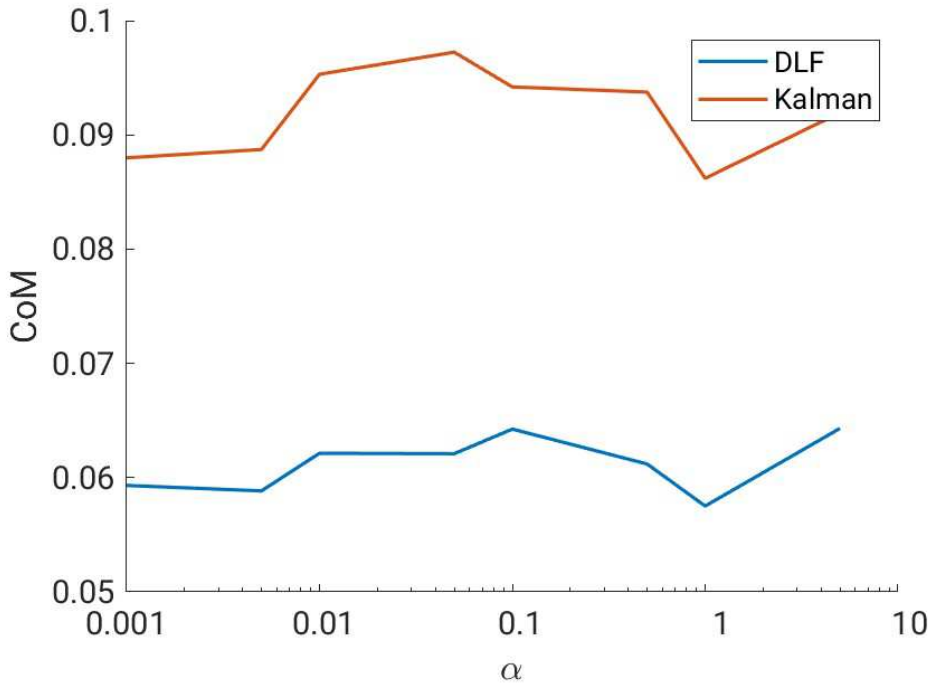

(c) CoM

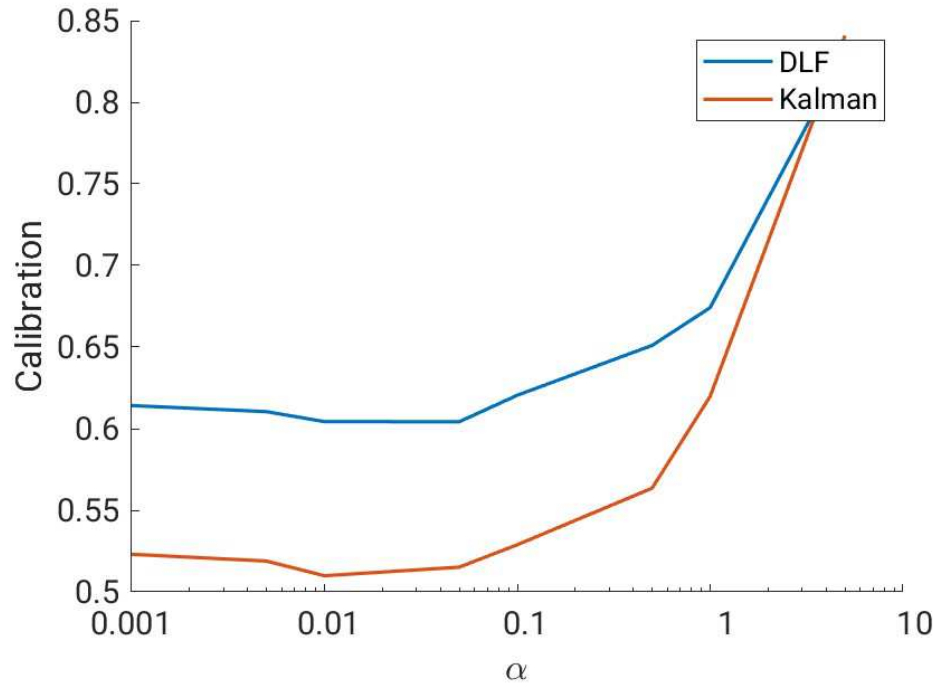

(d) Calibration

Fig. 16: Average (a) RMS, (b) Mass, (c) CoM errors and (d) Calibration of KF (blue), DLF (red), as a function of $\alpha$, across 50 runs for spatially correlated phase speed noise ($A = 0.05, \tilde{A} = 1$) and $I = 20$ randomly located observations at every observation time $T_O = \{0.05, 0.1, ..., 0.45\}$ with noisy initial data $\sigma \sim \mathcal{U}[\frac{1}{2}, \frac{3}{2}]$ and $\theta \sim \mathcal{U}[0, 1]$